\documentclass[a4paper,12pt]{amsart}

\usepackage{epsfig}
\usepackage{graphicx}

\usepackage[latin1]{inputenc}
\usepackage[T1]{fontenc}
\usepackage{indentfirst}
\usepackage{amssymb}
\usepackage{eufrak}
\usepackage{amsmath}
\usepackage{amsfonts}
\usepackage{amsthm}
\usepackage{mathrsfs}
\usepackage[bookmarks]{hyperref}
\usepackage{xypic}
\newcommand{\R}{\mathop{\rm Re}\nolimits}

\newcommand{\Log}{\mathop{\rm Log}\nolimits}
\newcommand{\val}{\mathop{\rm val}\nolimits}
\newcommand{\Val}{\mathop{\rm Val}\nolimits}
\newcommand{\ord}{\mathop{\rm ord}\nolimits}
\newcommand{\Arg}{\mathop{\rm Arg}\nolimits}
\newcommand{\Vol}{\mathop{\rm Vol}\nolimits}

\newcommand{\supp}{\mathop{\rm supp}\nolimits}
\newcommand{\Verte}{\mathop{\rm Vert}\nolimits}

\input epsf

\def \square{\smallskip \hfill \vrule width 5 pt height 7
pt depth - 2 pt \smallskip }

\newenvironment{prooof}
{\noindent {{\it Proof} \;}}{\hspace*{\fill}\square\vskip 8pt}

\theoremstyle{plain}
\newtheorem{Lem}[subsection]{Lemma}
\newtheorem{The}[subsection]{Theorem}
\newtheorem{Cor}[subsection]{Corollary}
\newtheorem{Pro}[subsection]{Proposition}
\theoremstyle{definition}
\newtheorem{Rem}[subsection]{Remark}

\newtheorem{Def}[subsection]{Definition}

\begin{document}

\title{Maximally Sparse Polynomials have Solid Amoebas}

\author{Mounir Nisse}
\address{Université Pierre et Marie Curie-Paris 6, IMJ (UMR 7586),
Labo: Analyse Algébrique, Office: 7C14, 175, rue du Chevaleret,\\ 75013
Paris, France}
\email{nisse@math.jussieu.fr}

\maketitle

\begin{abstract} Let  $f$ be an ordinary polynomial in
  $\mathbb{C}[z_1,\ldots ,z_n]$ with no negative exponents and with
  no factor of the form $z_1^{\alpha_1}\ldots z_n^{\alpha_n}$ where
  $\alpha_i$ are non zero natural integers. If we assume in addition
  that $f$ is a  maximally sparse polynomial (that its support is equal to the set
  of vertices of its Newton polytope), then a complement component of
  the amoeba $\mathscr{A}_f$ in $\mathbb{R}^n$ of the algebraic
  hypersurface $V_f\subset (\mathbb{C}^*)^n$ defined by $f$, has order
  lying in the support of $f$, which means that $\mathscr{A}_f$ is
  solid.
This gives an affirmative answer to Passare and Rullg\aa rd
  question in \cite{PR2-01}.

\end{abstract}

\setcounter{tocdepth}{1} \tableofcontents

\section{Introduction}

Mikael Passare and Hans Rullg\aa rd posed the following
question:

\vspace{0.3cm}

 "\, {\em Does every maximally sparse polynomial have a solid amoeba?}\, "

\vspace{0.3cm}

\noindent The purpose of this paper is to give an affirmative answer  to this
question. We use for this, Viro's patchworking principle applied
to the Passare and Rullg\aa rd function (see Section 2 for
definitions), Kapranov's theorem (see \cite{K-00}) and some properties of complex tropical hypersurfaces.
Note that here we can  assume that $f$ is a
polynomial with no negative exponent and with no factor of the form
$z^{\alpha}$  because our hypersurfaces lie in the algebraic
torus $(\mathbb{C}^*)^n$.

A polynomial $f$ is called {\em maximally sparse} if the support of $f$ is
equal to the set of the vertices of its Newton polytope $\Delta_f$
(see \cite{PR2-01}), in other word,  $f$ is a polynomial
with Newton polytope $\Delta$ with minimal number of monomials. An
amoeba $\mathscr{A}$ of degree $\Delta$ is called {\em solid }
if the number of connected component of $\mathbb{R}^n\setminus
\mathscr{A}$ is equal to the number of vertices of $\Delta$ , which is
the minimal number that an amoeba of degree $\Delta$ can have, (see \cite{PR1-04} or
\cite{R-01}). We prove the following theorem for any $n\geq 1$:

\begin{The} Let $V_f$ be an algebraic hypersurface in $(\mathbb{C}^*)^n$
  defined by a maximally sparse polynomial $f$. Then the amoeba
  $\mathscr{A}_f$ of $V_f$ is solid.
\end{The}

The paper is organized as follows. In section 2 we briefly review the definitions and the known results on tropical geometry and amoebas. We will then prove some properties of complex tropical hypersurfaces and we give a method for the construction of the set of arguments of a complex algebraic hypersurface defined by maximally sparse polynomial with Newton polytope a simplex in section 3. In section 4 we give the basic properties of Viro's local tropicalization. The proof of the main theorem will be given
 in section 5. It is based on  tropical localization of a special deformation of a complex structure on a hypersurface to the so-called by Grigory  Mikhalkin {\em complex tropical structure} which is the extrem possible degeneration.
In Appendix B  we give a geometric description of the set of arguments of the standard complex hyperplane, and finally in Appendix D we give an example which prove that maximally sparse polynomial is an optimal condition.
 
\vspace{0.3cm}

\noindent{\it Acknowledgment} The author would like to thank Professor Jean-Jacques Risler for his patient helps and helpful remarks, and Professor Mikael Passare for attracting my attention to the problem and
 useful discussions on the subject and others.

\section{Preliminaries}

In this paper we will consider only algebraic hypersurfaces $V$ in the
complex torus $(\mathbb{C}^*)^n$, where $\mathbb{C}^*
=\mathbb{C}\setminus \{ 0\}$ and $n\geq 1$ an integer. This means that $V$ is the zero locus of a Laurent polynomial:
$$
f(z) =\sum_{\alpha\in \supp (f)} a_{\alpha}z^{\alpha},\,\,\,\,\quad\quad\quad
z^{\alpha}=z_1^{\alpha_1}z_2^{\alpha_2}\ldots z_n^{\alpha_n}
\quad\quad\quad\quad\quad\quad\quad\quad  (1)
$$
where each $a_{\alpha}$ is a non-zero complex number and $\supp (f)$ is a
finite subset of $\mathbb{Z}^n$, called the support of the polynomial
$f$, with convex hull, in $\mathbb{R}^n$, the Newton polytope
$\Delta_f$ of $f$.
Moreover we assume that $\supp (f)\subset
\mathbb{N}^n$ and $f$ has no factor of the form $z^{\alpha}$ with
$\alpha =(\alpha_1,\ldots , \alpha_n)$.

 The {\em amoeba}  $\mathscr{A}_f$ of an algebraic hypersurface $V_f\subset
 (\mathbb{C}^*)^n$
 is by definition ( see M. Gelfand, M.M. Kapranov
 and A.V. Zelevinsky \cite{GKZ-94}) the image of $V_f$ under the map :
\[
\begin{array}{ccccl}
\Log&:&(\mathbb{C}^*)^n&\longrightarrow&\mathbb{R}^n\\
&&(z_1,\ldots ,z_n)&\longmapsto&(\log\mid z_1\mid,\ldots ,\log\mid
z_n\mid ).
\end{array}
\]

It was shown by M. Forsberg, M. Passare and A. Tsikh in \cite{FPT-00} that
there is an injective map between the set of components
$\{E_{\nu}\}$ of $\mathbb{R}^n\setminus \mathscr{A}_f$ and
$\mathbb{Z}^n\cap\Delta_f$:
$$
\ord :\{E_{\nu}\} \hookrightarrow \mathbb{Z}^n\cap\Delta_f
\qquad\qquad\qquad\qquad\qquad   (2)
$$

\begin{The}[Foresberg-Passare-Tsikh, (2000)] Each  component
  of $\mathbb{R}^n\setminus \mathscr{A}_f$  is a convex domain and there
  exists a locally constant function:
$$
\ord :\mathbb{R}^n\setminus \mathscr{A}_f \longrightarrow \mathbb{Z}^n\cap\Delta_f
$$
which maps different components of the complement of $\mathscr{A}_f$
to different lattice points of $\Delta_f$.
\end{The}

\vspace{0.3cm}

Let $\mathbb{K}$ be the field of the Puiseux series with real power, which is the field of the  series $\displaystyle{a(t) = \sum_{j\in A_a}\xi_jt^j}$ with $\xi_j\in \mathbb{C}^*$ and $A_a\subset \mathbb{R}$ is well-ordered set   with smallest element. It is well known that the field $\mathbb{K}$ is algebraically closed and of characteristic zero, and it has a non-Archimedean valuation $\val (a) = - \min A_a$:
\[ 
\left\{ \begin{array}{ccc}
\val (ab)&=& \val (a) + \val (b) \\
\val (a + b)& \leq& \max \{ \val (a)  ,\, \val (b)  \} ,
\end{array}
\right.
\]
and  we put $\val (0) = -\infty$. Let $f\in \mathbb{K}[z_1\ldots ,z_n]$ be a polynomial as in $(1)$ but the coefficients and the components of $z$ are in $\mathbb{K}$, and $V_{\mathbb{K}}$ be  the zero locus in $(\mathbb{K}^*)^n$ of the polynomial $f$. The following  piecewise affine linear function $\displaystyle{f_{trop} = \max_{\alpha\in\supp (f) } \{ \val (a_{\alpha}) + <\alpha , x>  \}}$ where $<,>$ is the
scalar product in $\mathbb{R}^n$
is called  a {\em tropical polynomial}.

\begin{Def} The tropical hypersurface $\Gamma_f$ defined by the tropical polynomial $f_{trop}$ is the subset of $\mathbb{R}^n$ image under  the valuation map of the algebraic hypersurface $V_{\mathbb{K}}$ over $\mathbb{K}$. 
\end{Def} 

\noindent We have the following Kapranov's theorem (see \cite{K-00}).

\begin{The}[Kapranov, (2000)] The tropical hypersurface $\Gamma_f$ is the set of points in $\mathbb{R}^n$ where the tropical polynomial  $f_{trop}$ is not smooth (called the corner locus of  $f_{trop}$).
\end{The}

\vspace{0.3cm}

{\bf Passare-Rullg\aa rd function.}

\vspace{0.3cm}

Let $A'$ be the subset of $\mathbb{Z}^n\cap\Delta_f$,  image of
$\{E_{\nu}\}$ under the order mapping $(2)$. M. Passare and
H. Rullg\aa rd proves in \cite{PR1-04} that the spine $\Gamma$ of the amoeba
$\mathscr{A}_f$ is given as a non-Archimedean amoeba defined by the
tropical polynomial
$$
f_{trop}(x) = \max_{\alpha\in A'}\{ c_{\alpha}+<\alpha , x>\} ,
$$
where $c_{\alpha}$ is defined by:
$$
c_{\alpha} = \R \left(   \frac{1}{(2\pi i)^n}\int_{\Log^{-1}(x)}\log
  \mid\frac{f(z)}{z^{\alpha}}\mid \frac{dz_1\wedge \ldots \wedge
    dz_n}{z_1\ldots z_n}\right)
\qquad\qquad\qquad  (3)
$$
where  $x\in E_{\alpha}$,\, $z = (z_1,\cdots ,z_n)\in (\mathbb{C}^*)^n$ and $<,>$ is the
scalar product in $\mathbb{R}^n$. In other words, the spine of
$\mathscr{A}_f$ is defined as the set of points in $\mathbb{R}^n$
where the piecewise affine linear function $f_{trop}$ is not
differentiable, or as the graph of this function where $\mathbb{R}$ is
the semi-field $(\mathbb{R}; \max , +)$.
Let us denote by $\tau$ the convex subdivision of $\Delta_f$ dual
to the tropical variety $\Gamma$.

\vspace{0.3cm}

\noindent
 We define the Passare-Rullg\aa rd's function on the Newton
polytope $\Delta_f$ as follows :

Let $\nu :\Delta_f\longrightarrow \mathbb{R}$ be the
function such that :
\begin{itemize}
\item[(i)]\, if $\alpha\in \Verte (\tau )$, then we set $\nu (\alpha
  )= -c_{\alpha}$
\item[(ii)]\, if $\alpha\in\Delta_v\setminus \Verte (\tau )$, where
  $\Delta_v$ is an element of the subdivision $\tau$ with maximal
  dimension, then we put $\nu (\alpha ) =<\alpha ,a_v>+b_v$, where
  $y=<x,a_v>+b_v$ is the equation of the hyperplane in
  $\mathbb{R}^n\times \mathbb{R}$ containing the points of coordinates
  $(\alpha , -c_{\alpha})\in\mathbb{R}^n\times \mathbb{R}$ for
  $\alpha\in \Verte (\Delta_v )$, \, $a_v=(a_{1,\, v},\ldots ,a_{
  n,\, v})\in \mathbb{R}^n$ and $b_v\in\mathbb{R}$.
\end{itemize}

If $f$ is the polynomial given by $(1)$, we define a family
of polynomials $\{ f_t\}_{t\in ]0,\frac{1}{e}]}$ as follows :
$$
f_t(z) = \sum_{\alpha\in \supp (f)} \xi_{\alpha}t^{\nu
(\alpha )}z^{\alpha}\qquad\qquad\qquad\qquad\qquad (4)
$$
where $\xi_{\alpha} = a_{\alpha} e^{\nu (\alpha )}$.

\section{Complex tropical hypersurfaces.}

 Let $h$ be  a strictly positive real number and $H_h$ be the self
 diffeomorphism of $(\mathbb{C}^*)^n$ defined by :
\[
\begin{array}{ccccl}
H_h&:&(\mathbb{C}^*)^n&\longrightarrow&(\mathbb{C}^*)^n\\
&&(z_1,\ldots ,z_n)&\longmapsto&(\mid z_1\mid^h\frac{z_1}{\mid
  z_1\mid},\ldots ,\mid z_n\mid^h\frac{z_n}{\mid z_n\mid} ),
\end{array}
\]
which defines a new complex structure on $(\mathbb{C}^*)^n$
denoted by $J_h = (dH_h)^{-1}\circ J\circ (dH_h)$ where $J$ is the
standard complex structure.

\noindent A $J_h$-holomorphic hypersurface $V_h$ is a hypersurface
holomorphic with respect to the $J_h$ complex structure on
$(\mathbb{C}^*)^n$. It is equivalent to say that $V_h = H_h(V)$ where
$V\subset (\mathbb{C}^*)^n$ is an holomorphic hypersurface for the
standard complex structure $J$ on $(\mathbb{C}^*)^n$.

Recall that the Hausdorff distance between two closed subsets $A,
B$ of a metric space $(E, d)$ is defined by:
$$
d_{\mathcal{H}}(A,B) = \max \{  \sup_{a\in A}d(a,B),\sup_{b\in B}d(A,b)\}.
$$
Here we take $E =\mathbb{R}^n\times (S^1)^n$, with the distance
defined as the product of the
Euclidean metric on $\mathbb{R}^n$ and the flat metric on $(S^1)^n$.

\begin{Def} A complex tropical hypersurface $V_{\infty}\subset
  (\mathbb{C}^*)^n$ is the limit (with respect to the Hausdorff
  metric on compact sets in $(\mathbb{C}^*)^n$) of a  sequence of a
  $J_h$-holomorphic hypersurfaces $V_h\subset (\mathbb{C}^*)^n$ when
  $h$ tends to zero.
\end{Def}

Using Kapranov's theorem \cite{K-00}, Mikhalkin gives an algebraic
definition of a complex tropical hypersurfaces (see \cite{M2-04}) as
follows :\\

If $ a \in \mathbb{K}^*$ is the Puiseux series
$\displaystyle{a = \sum_{j\in A_a}\xi_jt^j}$ with $\xi\in \mathbb{C}^*$
and $A_a\subset \mathbb{R}$ is a well-ordered set with smallest element,
then we have a non-Archimedean valuation on $\mathbb{K}$ defined by
$\val (a ) = - \min A_a$. We complexify the valuation
map as follows :
\[
\begin{array}{ccccl}
w&:&\mathbb{K}^*&\longrightarrow&\mathbb{C}^*\\
&&a&\longmapsto&w(a ) = e^{\val (a )+i\arg (\xi_{-\val
    (a )})}
\end{array}
\]
Let  $\Arg$ be the argument map $\mathbb{K}^*\rightarrow S^1$ defined by: for any $ a \in \mathbb{K}$ a Puiseux series so that
$\displaystyle{a = \sum_{j\in A_a}\xi_jt^j}$, then $\Arg (a) = e^{i\arg (\xi_{-\val (a)})}$ (this map extends the map $\mathbb{C}^*\rightarrow S^1$ defined by $\rho e^{i\theta} \mapsto e^{i\theta}$ which we denote by $\Arg$).

Applying this map coordinate-wise we obtain a map :
\[
\begin{array}{ccccl}
W:&(\mathbb{K}^*)^n&\longrightarrow&(\mathbb{C}^*)^n
\end{array}
\]

\begin{The}[Mikhalkin, (2002)]  The set  $V_{\infty}\subset (\mathbb{C}^*)^n$
  is a complex tropical hypersurface if and only if there
  exists an algebraic hypersurface
  $V_{\mathbb{K}}\subset(\mathbb{K}^*)^n$ over $\mathbb{K}$ such that
  $\overline{W(V_{\mathbb{K}})} = V_{\infty}$, where $\overline{W(V_{\mathbb{K}})}$ is the closure of $W(V_{\mathbb{K}})$ in  $(\mathbb{C}^*)^n \approx \mathbb{R}^n\times (S^1)^n$ as  a Riemannian manifold with metric defined by the standard Euclidean metric of $\mathbb{R}^n$ and  the standard  flat metric of the torus.
\end{The}

Recall that we have the following commutative diagram:
\begin{equation}
\xymatrix{
(\mathbb{K}^*)^n\ar[rr]^{W}\ar[dr]_{\Log_{\mathbb{K}}}&&(\mathbb{C}^*)^n\ar[dl]^{\Log}\cr
&\mathbb{R}^n
}\nonumber
\end{equation}
where $\Log_{\mathbb{K}}(z_1,\ldots ,z_n)=(\val (z_1),\ldots ,
\val (z_n))$, which means that $\mathbb{K}$ is equipped with the
norm defined by $\mid z\mid_{\mathbb{K}}= e^{\val (z)}$ for any
$z\in\mathbb{K}^*$.

Let $V_{\infty}\subset(\mathbb{C}^*)^n$ be a complex tropical
hypersurface of degree $\Delta$. This means that
$ V_{\infty} =
\overline{W(V_{\mathbb{K}})}$ where $V_{\mathbb{K}}\subset (\mathbb{K}^*)^n$ is
an algebraic hypersurface over $\mathbb{K}$ defined by the
non-Archimedean polynomial $\displaystyle{f_{\mathbb{K}}(z) = \sum
  _{\alpha\in \Delta\cap \mathbb{Z}^n}a_{\alpha}z^{\alpha}}$. By Kapranov's theorem (see \cite{K-00}),
$\Gamma = \Log_{\mathbb{K}}(V_{\mathbb{K}})$ is a tropical hypersurface (called non-Archimedean amoeba associated to the polynomial $f_{\mathbb{K}}$ and denoted  by $\mathscr{A}_{f_{\mathbb{K}}}$); we denote by $\tau$ the
subdivision of $\Delta$ dual to $\Gamma$.

\begin{Def} The  complex numbers $w(a_{\alpha})$ are called the  complex
  tropical coefficients "defined" by $V_{\infty}$. They are well
  defined if we suppose that for some fixed index $\alpha_0\in\Delta$,
  $w(a_{\alpha_0})=1$
\end{Def}

In general we have the following (see Mikhalkin \cite{M2-04} for
$n=2$):

\begin{Pro}Let $V_{\mathbb{K}}\subset (\mathbb{K}^*)^n$ as above. Then
for any  two indices  $\alpha$ and $\beta$ in $\Verte (\tau
  )$, the quotients $\frac{w(a_{\alpha})}{w(a_{\beta})}$ are well
  defined and depend only on $W(V_{\mathbb{K}})$.
\end{Pro}

\begin{prooof} First of all , we may assume that $\alpha$ and $\beta$ are
  adjacent to the same edge $E$ of $\tau$ and we proceed by
  induction on vertices. Secondly, using an automorphism of $(\mathbb{C}^*)^n$ if
  necessary, we may assume  that $E=[0,k]\times\{
  0\}\subset\mathbb{R}^n$. Let $E^*\subset \Gamma$ be the dual of $E$ and
  $U\subset \mathbb{R}^n$ a
small
 neighborhood of a point $x\in Int
  (E^*)$,
then we have
  $\Log^{-1}(U)\cap V_{\infty} = \Log^{-1}(U)\cap V_{\infty ,\, E}$
  where $V_{\infty ,\, E}$ is the complex tropical hypersurface
  defined in the same way of $V_{\infty}$ but by taking the truncation
  of $f$ to $E$. Indeed,
  the tropical monomials corresponding to lattice points in $E$
  dominate the tropical monomials corresponding to lattice points in
  $\Delta\setminus E$ (it's Kapranov's theorem \cite{K-00}).
  Hence we can assume that $\Delta = E$ and
  prove the result for $E$.

\noindent Let $f_{\mathbb{K}}^E(z_1) = a_0z_1^k +a_1z_1^{k-1}+\ldots
  +a_k\in \mathbb{K}[z_1]$ be a non-Archimedean polynomial in one
  variable such that $W(V_{f_{\mathbb{K}}^E}) = V_{\infty ,\, E}$ (it
  can be seen as the truncation of $f_{\mathbb{K}}$ on $E$). The field
  $\mathbb{K}$ is algebraically closed, hence the polynomial
  $f_{\mathbb{K}}^E$ has $k$ roots $r_1,\ldots ,r_k$ in $\mathbb{K}$
  such that $\displaystyle{\prod_{j=1}^kr_j = (-1)^k\frac{a_k}{a_0}}$. On
  the other hand $V_{\infty ,\, E}$ is the union of
 subsets
  $\displaystyle{\cup_{i=1}^s\mathscr{C}_i}$ in  $(\mathbb{C}^*)^n$
 defined by $\mathscr{C}_i = \{ (z_1,\ldots ,z_n)\in (\mathbb{C}^*)^n/\,\,
 \prod_{j=1}^{k_i} (z_1 - c_{ij})= 0 \}$ 
such that for any $i$ we have $\log\mid
  c_{ij}\mid = c_i$  where $c_i$ are  a constants depending only on $V_{\infty ,\, E}$.
 Indeed,  $\Log (V_{\infty ,\, E})$ is an
  hyperplane in $\mathbb{R}^n$ orthogonal to the $x_1$-axis, and $k=\sum
  k_i$. Any $w(r_j)\in \mathscr{C}_i$ for some $i$, so there exists $j$ such that
  $w(r_j)=c_{ij}$. This means that $w(r_j)$ is a solution of
  the equation $z_1-c_{ij}=0$ (in the field of the complex
  numbers).

 Then
  $\displaystyle{\prod_{j=1}^{k}w(r_j) =\prod_i \prod_{j=1}^{k_i}c_{ij}}$ and
  hence we have:

$$
\frac{w(a_k)}{w(a_0)} = (-1)^k \prod_i \prod_{j=1}^{k_i}c_{ij}
$$
which depends only on $V_{\infty ,\, E}$ and hence only on $V_{\infty}$; this
  proves Proposition 3.4.


\end{prooof}

Let $f_{\mathbb{K}}$ be a  polynomial in
$\mathbb{K}[z_1,\ldots ,z_n]$ with Newton polytope a simplex $\Delta$
such that $\supp (f_{\mathbb{K}}) = \Verte (\Delta )$;  this implies that the  corresponding
non-Archimedean amoeba $\mathscr{A}_{f_{\mathbb{K}}}$
   has only one vertex. Assume that there exists $\{ g_{\mathbb{K},\,
  u}\}_{u\in [0,1]}$ a family of non-Archimedean polynomials
defined by
$\displaystyle{ g_{\mathbb{K},\, u}(z) =
    f_{\mathbb{K}}(z) + \sum_{\beta\in A}a_{\beta ,\, u}z^{\beta}}$
 where $A\subset (\Delta \cap\mathbb{Z}^n)\setminus \Verte (\Delta
    )$
satisfy the following properties :
\begin{itemize}
\item[(i)]\, the complement components of  the non-Archimedean amoeba
  $\mathscr{A}_{g_{\mathbb{K},\, 1}}$ of
  $g_{\mathbb{K},\, 1}$ are in bijection with $\Verte (\Delta )\cup A$
  by the order map and if  we denote by $\tau$ the  subdivision of $\Delta$
dual to the non-Archimedean amoeba $\mathscr{A}_{g_{\mathbb{K},\, 1}}$
  we assume that $\tau$ is a triangulation,
\item[(ii)]\, let $\nu$ be the Passare-Rullg\aa rd
  function associated to the amoeba of $f_{\mathbb{K}}$ (in this case $\nu (\alpha ) = -\log\mid a_{\alpha}\mid$ for any $\alpha\in \Verte (\Delta )$), and 
for any $\beta\in A$ and $0\leq u\leq 1$, $\val
  (a_{\beta ,\, u}) = (1-u)\val
  (a_{\beta ,\, 0})  + u\val
  (a_{\beta ,\, 1})$, where $\val
  (a_{\beta ,\, 0}) =   -\nu (\beta
  )$,
 and   $\arg
  (a_{\beta ,\, u}) = \arg (a_{\beta ,\, 1})$ for each $u$ such that $0\leq u\leq 1$.
\end{itemize}

\vspace{0.2cm}

Let us denote by $D_{std}$ the lift set in $\mathbb{R}^n$ of the argument of the complex tropical
hyperplane $W(H)$ where $H$
is the hyperplane in $(\mathbb{K}^*)^n$
 defined
by the polynomial $z_1+\cdots +z_n+1 = 0$,
with degree the standard simplex $\Delta_{std} = \{
(x_1,\ldots ,x_n)\in\mathbb{R}^n \, \mid \, x_j\geq 0,\, x_1+\cdots
+x_n\leq 1\}$.

\begin{Pro}Let $f_{\mathbb{K}}$ and $g_{\mathbb{K}
  ,\, u}$ having the above properties. Then
$W (V_{f_{\mathbb{K}}}) =  W (V_{g_{\mathbb{K},\,
      0}})$ if and only if $A$ is empty.
\end{Pro}

We can remark that if $A$ is empty, $W (V_{f_{\mathbb{K}}}) =  W (V_{g_{\mathbb{K},\,
      0}})$   because $
g_{\mathbb{K},\, u} = f_{\mathbb{K}} $ in this case. Let $\Delta_i$, for $1\leq i\leq s$,
be the simplices of dimension $n$ of the triangulation $\tau$. We have the following diagram:

\begin{equation}
\xymatrix{
&&\mathbb{R}^n\ar[d]^{\pi}\cr
(\mathbb{K}^*)^n\ar[r]^{W}&(\mathbb{C}^*)^n\ar[d]^{\Log}\ar[r]^{\Arg}&(S^1)^n\cr
&\mathbb{R}^n&
}
\end{equation}
where the map $\pi : \mathbb{R}^n\longrightarrow (S^1)^n $ is the projection of the universal covering of the torus, and $\Arg (\rho e^{i\theta}) = e^{i\theta}$: see above.

\begin{Lem} Let $f_{\mathbb{K}}$ and $g_{\mathbb{K}
  ,\, u}$ with  properties (i) and (ii). Then
there exist invertible matrices $\displaystyle{\{
  L_i\}_{i=1}^{s}\subset GL(n,\mathbb{R})}$ with
  coefficients in $\mathbb{Z}$ and positive determinant depending
  only on the triangulation $\tau$ (where $s$ is the number of element
  of $\tau$), and real
  vectors $\displaystyle{\{ (v_i)_{i=1,\ldots ,s}\mid  v_i\in \mathbb{R}^n\}}$  depending only on  the
  complex tropical hypersurface $W ( V_{g_{\mathbb{K},\, 0}})$
 such that :\\
 if $v$ is the only vertex of the non-Archimedean
  amoeba $\mathscr{A}_{f_{\mathbb{K}}}$ we have\\
$\displaystyle{\Arg (\Log^{-1}(v)\cap W(V_{g_{\mathbb{K},\, 0}})) =
  \cup_{i=1}^{s}\mathscr{C}_i}$ where $\mathscr{C}_i = (tr_{v_i}\circ
  {}^tL_i^{-1}(D_{std})\setminus \mathcal{R}_i)/(2\pi \mathbb{Z})^n$ and $tr_{v_i}$
  are translations,  with $\mathcal{R}_i/(2\pi \mathbb{Z})^n\subset \Arg (W(V_{g_{\mathbb{K},\, 1}^{\Delta_i}}))$
and  depends only on the coefficients of $g_{\mathbb{K},\, 1}$ with index in $\Delta_j$'s
which  has a common face with $\Delta_i$.
\end{Lem}

\begin{prooof} We do not need the case $n=1$ because Theorem 1.1 is obviously true for $n=1$. However, we postpone the proof of Lemma 3.6. and Proposition 3.5. for $n=1$ in the Appendix.

\noindent\,\, {\it Case} $n\geq 2$.

Let $f_{\mathbb{K}}(z)=a_0+\sum_{j=1}^na_jz_1^{\alpha_{1j}}\ldots
z_n^{\alpha_{nj}}$  
and $V_{f_{\mathbb{K}}}$ its
zero locus in $(\mathbb{K}^*)^n$. So $V_{f_{\mathbb{K}}}$ is the image
of a hyperplane in $(\mathbb{K}^*)^n$ by the endomorphism
 $\tilde{L}_{\Delta} : (\mathbb{K}^*)^n\longrightarrow
(\mathbb{K}^*)^n$ defined by the change of the variable $(z_1',\ldots ,z_n') \leadsto
(z_1,\ldots ,z_n)$,  $z_j'= z_1^{\alpha_{1j}}\ldots
z_n^{\alpha_{nj}}$
for $j=1,\ldots ,n$.

This gives an
endomorphism of the rings:
\[
\begin{array}{ccc}
\mathbb{K}[ {z'}_1^{\pm 1},\ldots , {z'}_n^{\pm 1}]&\longrightarrow&\mathbb{K}[ {z_1}^{\pm 1},\ldots , {z_n}^{\pm 1}]\\
z_j'&\longmapsto&z_1^{\alpha_{1j}}\ldots
z_n^{\alpha_{nj}}
\end{array}
\]
Let ${}^tL_{\Delta}$ be  the transpose matrix of the linear part
$L_{\Delta}$ of the affine linear surjection
which transform the standard simplex (i.e. with the $n+1$ vertices $(0,\ldots ,0),\,
(1,0,\ldots ,0)$,\,\, $(0,\ldots ,0,1,0,\ldots ,0)$ and $(0,\ldots
,0,1)$) to $\Delta$. 
So we  obtain the automorphism of $\mathbb{R}^n$ defined by:
\[
\begin{array}{ccccl}
^{t}L_{\Delta}&:&\mathbb{R}^n&\longrightarrow&\mathbb{R}^n\\
&&\left( \begin{array}{c}
\val (z_1)\\ \vdots\\\val (z_n)
\end{array}\right)&\longmapsto&
\left( \begin{array}{ccc}
\alpha_{11}&\hdots&\alpha_{n1}\\
\vdots&\ddots&\vdots\\
\alpha_{1n}&\hdots&\alpha_{nn}
\end{array}\right)
\left( \begin{array}{c}
\val (z_1)\\ \vdots\\\val (z_n)
\end{array}\right)  =
\left( \begin{array}{c}
\val (z_1')\\\vdots\\ \val (z_n')
\end{array}\right)
\end{array}
\]
and we have an homomorphism $\tilde{L}_{\Delta}$ of the multiplicative group
$(\mathbb{K}^*)^n$ such that
 the following diagram is commutative:
\begin{equation}
\xymatrix{
(\mathbb{K}^*)^n\ar[d]_{\Val}&&(\mathbb{K}^*)^n\ar[ll]_{\tilde{L}_{\Delta}}\ar[d]^{\Val}\cr
\mathbb{R}^n\ar[rr]^{{}^tL_{\Delta}^{-1}}&&\mathbb{R}^n.
}\nonumber
\end{equation}
 Then the image of the hyperplane in  $(\mathbb{K}^*)^n$  defined by the
polynomial $f_{\mathbb{K},\, ltd}(z')=a_0+\sum_{j=1}^na_j{z'}_j$ is precisely the
hypersurface defined by $f_{\mathbb{K}}$.
Here we have
$$
L_{\Delta} = \left(
\begin{array}{ccc}
\alpha_{11}&\hdots&\alpha_{1n}\\
\vdots&\ddots&\vdots\\
\alpha_{n1}&\hdots&\alpha_{nn}
\end{array}
\right) .
$$
So we obtain the following commutative diagram:

\begin{equation}
\xymatrix{
(\mathbb{K}^*)^n\ar[d]_{\tilde{W}}&&(\mathbb{K}^*)^n\ar[d]^{\tilde{W}}\ar[ll]_{\tilde{L}_{\Delta}}\cr
\mathbb{C}^n\ar[rr]^{\overline{L}_{\Delta}}\ar[d]_{\exp}&&\mathbb{C}^n\ar[d]^{\exp}\cr
(\mathbb{C}^*)^n\ar[rr]^{M_{\Delta}}&&(\mathbb{C}^*)^n
}\nonumber
\end{equation}
where $\tilde{W}= (\Val ,\Arg )$ (i.e. $\tilde{W}(z)= (\val (z) + i\arg (z))$), and $M_{\Delta}$ is the
endomorphism of $(\mathbb{C}^*)^n$ covered by
$\overline{L}_{\Delta} = ({}^tL_{\Delta})^{-1}\otimes \mathbb{C}$. Using the map $W= 
\exp\circ \tilde{W}$, we obtain
the following commutative diagram:
\begin{equation}
\xymatrix{
(\mathbb{K}^*)^n\ar[d]_{W}&&(\mathbb{K}^*)^n\ar[d]^{W}\ar[ll]_{\tilde{L}_{\Delta}}\cr
(\mathbb{C}^*)^n\ar[rr]^{M_{\Delta}}\ar[d]_{\Log}&& (\mathbb{C}^*)^n\ar[d]^{\Log}\cr
\mathbb{R}^n\ar[rr]^{{}^tL_{\Delta}^{-1}}&&\mathbb{R}^n
}\nonumber
\end{equation}
The degree of the map $M_{\Delta}$ is equal to the determinant of
$L_{\Delta}$ i.e. $deg (M_{\Delta}) = \det (L_{\Delta})= n!\Vol
(\Delta )$. Let $H_a$ be the  hyperplane in $(\mathbb{K}^*)^n$ defined by the
polynomial $f_{\mathbb{K},\,
    ltd}(z') =
a_0+\sum_{j=1}^na_j{z'}_j$,
then  $W(V_{f_{\mathbb{K}}})\, =\, M_{\Delta}(W(H_a))$.

\vspace{0.2cm}

{\it Claim} {\bf 1}: Let $f_{\mathbb{K},\,
    std}(z) = 1+\sum_{j=1}^nz_j$. Then we have $W(H_a) = \tau_v\circ W(H)$ where $\tau_v$ denotes the translation in the
multiplicative algebraic torus $(\mathbb{C}^*)^n$ by an element $v\in
(\mathbb{C}^*)^n$ well defined by the coefficients of $f_{\mathbb{K},\,
    ltd}$.

\begin{prooof} We are in the algebraic torus, then we can assume that
$a_0=1$ and the valuation of each other coefficients is zero. Indeed,
let $\Phi_a$ be the automorphism of $(\mathbb{K}^*)^n$ defined by $\Phi_a
(z_1,\ldots ,z_n) = ( t^{\val (a_1)}z_1,\ldots ,t^{\val (a_n)}z_n)$,
then $f_{\mathbb{K},\, ltd}\circ \Phi_a$ has the required assertion
and $\Arg (W(H_a)) = \Arg (W(V_{f_{\mathbb{K},\,
    ltd}\circ\Phi_a}))$. We can see that if $f_{\mathbb{K},\,
  ltd}\circ \Phi_a(z) = 1+\sum_{i=0}^na_i'z_i$ then $\Arg (W(V_{f_{\mathbb{K},\,
    ltd}\circ\Phi_a})) = \tau_v (\Arg (W(H)))$ where
$\tau_v$ is the multiplication in the real  torus
$(S^1)^n$ by $(e^{i\arg (a_1')},\ldots ,e^{i\arg (a_n')})$.
\end{prooof}

\noindent Let $\tilde{v} = {}^tL_{\Delta}^{-1}(v)$ (where $v$ is viewed as a vector in the universal covering), then  we have
$W(V_{f_{\mathbb{K}}})\, =\, M_{\Delta}(W(H_a)) =
 \tau_{\tilde{v}}\circ M_{\Delta}(W(H))$; so we obtain
 $\Arg (W(V_{f_{\mathbb{K}}})) = (\tau_{\tilde{v}}\circ
{}^tL_{\Delta}^{-1}(D_{std}))/(2\pi\mathbb{Z})^n$
where $D_{std}\subset \mathbb{R}^n$ is the lift set of the argument
of  $W(H)$ in the universal covering of
$(S^1)^n$ and $\tilde{v} = {}^tL_{\Delta}^{-1}(\arg (a_1'),\ldots ,\arg (a_n'))\in
\mathbb{R}^n$.

\vspace{0.2cm}

\noindent We can remark that  for any $\Delta_i\in \tau$,
the argument of $W(V_{g_{\mathbb{K},\,
    u}^{\Delta_i}})$
is independent of $u$,
because the deformation is given such
that the combinatorial type of the tropical hyperplane $\Log\circ W(V_{g_{\mathbb{K},\,
    u}^{\Delta_i}})$ is the same for any $u$ and the argument  of the
coefficients of $g_{\mathbb{K},\, u}^{\Delta_i}$ are independent of $u$
by construction.
We denote by $ad (i)$ the set of $j$ so that $\Delta_j$ is adjacent to $\Delta_i$.
 Let $\mathcal{R}_{ij}$ be the subset of the lift in $\mathbb{R}^n$ of
$ \Arg (W(V_{g_{\mathbb{K},\, 1}^{\Delta_i}})) \cap  \Arg (W(V_{g_{\mathbb{K},\,
    1}^{\Delta_j}}))$ not in the lift of  $\Arg (W(V_{g_{\mathbb{K},\,
    1}^{\Delta_i \cup \Delta_j}}))$, and put $\displaystyle{\mathcal{R}_i = \cup_{j\in ad (i)} \mathcal{R}_{ij}}$.
The $\mathcal{R}_{ij}$ depends only on the coefficients of $g_{\mathbb{K},\, 1}$ with index in $\Delta_i\cup \Delta_j$.
Hence, at the limit (i.e. $u=0$), we have:

\begin{eqnarray}
\Arg (W(V_{g_{\mathbb{K},\, 0}})) &=&
\bigcup_{\Delta_i\in\tau} \{ \Arg ((W(V_{g_{\mathbb{K},\,
    1}^{\Delta_i}}))\setminus \mathcal{R}_i)/(2\pi \mathbb{Z})^n\}\nonumber \\
    &=& \bigcup_{\Delta_i\in\tau} (\tau_{\tilde{v}_i}\circ
{}^tL_{\Delta_i}^{-1}(D_{std})\setminus \mathcal{R}_i)/(2\pi\mathbb{Z})^n \nonumber
\end{eqnarray}
where $L_{\Delta_i}$ is the linear part of the affine linear surjection map
between the standard simplex and $\Delta_i$, $\tilde{v}_i$ are the translation vectors
on the torus corresponding to the truncations $g_{\mathbb{K},\,
    1}^{\Delta_i}$
as described above.
\end{prooof}

\noindent {\it Proof of  Proposition 3.5.} For each $u$, let us take the following notations :\,  $V_{\infty ,\, u} = W(V_{g_{\mathbb{K},\, u}})$ and $V_{\infty ,\, f} = W(V_{f_{\mathbb{K}}})$.
The principal arguments of Proposition 3.5,  are the fact that, firstly the non-Archimedean amoeba $\Gamma_{\infty}$ has only one vertex. Secondly, if $f_{\mathbb{K}}$ is maximally sparse, then the lifting of the boundary $\partial \overline{\Arg (V_{\infty ,\, f})}$ of the closure of the set of argument  of the complex tropical hypersurface $V_{\infty ,\, f}$ (called by M. Passare the coamoeba of the complex tropical hypersurface $V_{\infty ,\, f}$  and denoted by $co\mathscr{A}_{V_{\infty ,\, f}}$),  are the hyperplanes orthogonal to the edges $E_{\alpha_i \alpha_j}$ of $\Delta$; in addition 
$\overline{\Arg (V_{\infty ,\, 0})}$ contains extra-pieces,
 where $V_{\infty ,\, 0} = \lim_{u\rightarrow 0}V_{\infty ,\, u}$.
Indeed,  let $H_{ij}$ be the hyperplane  image under ${}^tL^{-1}$ of the hyperplane $H_{ij}^{std}$ orthogonal to the edge $E_{\alpha_i \alpha_j}^{std}$ of the standard simplex such that $E_{\alpha_i \alpha_j} = L (E_{\alpha_i \alpha_j}^{std})$. Then we have $<{}^tL_{\Delta}^{-1}(H_{ij}^{std}), E_{\alpha_i \alpha_j}> = <H_{ij}^{std}, L^{-1}(E_{\alpha_i \alpha_j})> = <H_{ij}^{std}, E_{\alpha_i \alpha_j}^{std}> =0$. Secondly,
each edge $E_{\alpha_i\alpha_j}$ of $\Delta$ is dual to an $(n-1)$-polyhedron $E^*_{ij}\subset \Gamma_{\infty} = \mathscr{A}_{f_{\mathbb{K}}}\subset \mathbb{R}^n$. Let $x\in E^*_{ij}$ and $U$ be a small ball in $\mathbb{R}^n$ centered at $x$. Then, using Kapranov's theorem \cite{K-00}, we have:
$$
\overline{\Arg (\Log^{-1}(U)\cap V_{\infty ,\, f})} \subset \mathcal{N}_{\varepsilon} ( \Arg \{ z\in (\mathbb{C}^*)^n / \, a_{\alpha_i}z^{\alpha_i} + a_{\alpha_j}z^{\alpha_j} = 0\}),
$$
where $\alpha_i$ and $\alpha_j$ are the vertices of $\Delta$ bounding the edge $E_{\alpha_i \alpha_j}$ and $\mathcal{N}_{\varepsilon}$ designate the $\varepsilon$-neighborhood . Hence we obtain:
$$
\arg (a_{\alpha_i}) + <\alpha_i,\Arg (z)> = \pi + \arg (a_{\alpha_j}) + <\alpha_j,\Arg (z)> + 2k\pi ,
$$
where $k\in\mathbb{Z},\,$ $\Arg (z) = (\arg (z_1),\ldots ,\arg (z_n)),\, \alpha_i = (\alpha_{i1},\ldots , \alpha_{in})$,\, $\alpha_j = (\alpha_{j1},\ldots , \alpha_{jn})$, and $<,>$ is the Euclidean scalar product. So the hyperplanes $H_{ij}$ in $\mathbb{R}^n$ of equations:
$$
\arg (a_{\alpha_i}) -\arg (a_{\alpha_j}) + \sum_{l=1}^n(\alpha_{il} - \alpha_{jl}).x_l = (2k+1)\pi,
$$
are the boundary of the set $\overline{\Arg (V_{\infty ,\, f})}$ because the set of arguments of $V_{\infty ,\, f}$ can be only in one side of $H_{ij}$.
Let us describe now the boundary of $\overline{\Arg (V_{\infty ,\, 0})}$.

\begin{Lem} For any $u\geq 0$, there are extra-pieces $\mathcal{P}_{j,\, u}$ contained in $\Arg (V_{\infty ,\, u})$ with no vanishing volume, such that $\mathcal{P}_{j,\, u} \cap \Arg (V_{\infty ,\, f}) = \phi$.
\end{Lem}

\begin{prooof} Let us call the hyperplanes $H_{ij}$ {\em external hyperplane}, and assume that $A$ contains just one point $\beta$ in the interior of $\Delta$. Let
 $V_{\infty ,\,g^{\Delta_i}_u}$ for $i=0,\ldots , n-1$ be the complex tropical hypersurface image under the map $W$
of the hypersurface in $(\mathbb{K}^*)^n$ defined by the truncation of $g_{\mathbb{K},\, u}$ to
$\Delta_{i}$ where  $\Delta_i$ is an element of the triangulation $\tau =\{ \Delta_0,\ldots ,\Delta_n\}$. There exists an external hyperplane $H_{rs}\subset \partial \overline{\Arg (V_{\infty ,\, f})}$ intersecting the union $\displaystyle{\cup_{i=0}^n\Arg (V_{\infty ,\,g^{\Delta_i}_u})}$ in its interior. Because otherwise, this means that each connected component of the complement of $\Arg (V_{\infty ,\, f})$ in the torus $(S^1)^n$ is strictly contained in some complement component of $\Arg (V_{\infty ,\, g^{\Delta_i}_u })$ (indeed, if that inclusion is not strict, then there is at least one face of some  complement component of $\Arg (V_{\infty ,\, f})$ intersecting the interior of $\Arg (V_{\infty ,\, g^{\Delta_i}_u })$), and then  $\Vol ((S^1)^n\setminus \Arg (V_{\infty ,\, g^{\Delta_i}_u })) > \Vol ((S^1)^n\setminus \Arg (V_{\infty ,\, f}))$; that  contradicts the fact that the volume of the last two sets is the same (see Appendix B).

%
%
%
%

 Hence
$\Arg (V_{\infty ,\, u})$ contains some extra-pieces $\mathcal{P}_{j,\, u}$ in the exterior of $\Arg (V_{\infty ,\, f})$, because the set of argument of $V_{\infty ,\, f}$ can be only in one side of the external hyperplanes. Let $\displaystyle{\mathcal{P}_j = \lim_{u\rightarrow 0}\mathcal{P}_{j,\, u}}$, then $\Vol (\mathcal{P}_j)\ne 0$. Indeed, if $\Vol (\mathcal{P}_j) = 0$, this means that the valuation of the coefficient with index $\beta$ tends to  $-\infty$ (in other word, this means that the coefficient with index $\beta$ tends to zero), which is not the case by construction, because the valuation of that coefficient tends to $-\nu_{PR}(\beta )$, which is finite. 
\end{prooof}

On the other hand we have: $\displaystyle{\Arg (V_{\infty ,\, 0}) = \lim_{u\rightarrow 0} \Arg (V_{\infty ,\, u})}$. So the boundary of $\overline{\Arg (V_{\infty ,\, 0})}$ contains other pieces in $\partial \mathcal{P}_j$, not in the hyperplanes $H_{ij}$. Hence $\Arg (V_{\infty ,\, f})$ cannot be equal to $\Arg (V_{\infty ,\, 0})$.

\noindent If $A$ contains  more than one point, we use induction on the cardinality of $A$, and we subdivide $\Delta$ into  at most $n+1$ simplex with common vertex $\beta\in A$. Using the same reasoning  we have the result.

\begin{Rem}
 For any $u$, there are extra-pieces $\mathcal{P}_{j,\, u}$ with no vanishing volume in $\Arg (V_{\infty ,\, u})$ (see for example figure 2 and 3 for $n=2$),  corresponding to the dual of the edges of the subdivision of $\Delta$ (dual to $\Gamma_{\infty ,\, u}$) other than the edges of $\Delta$. So the sets $\Arg
(V_{\infty ,\, 0})$ and $\Arg
(V_{\infty ,\, f})$ cannot be equal even if $u$ tends to some negative real number,  this means even if the valuation of the coefficients $a_{\beta}$'s are above the hyperplane in $\mathbb{R}^{n+1}$ passing through the points of coordinates $(\alpha ,\nu (\alpha ))$ with $\alpha$ in $\Verte (\Delta )$. If we add in the hypothesis of the Proposition 3.5  that the 
non-Archimedean amoeba $\mathscr{A}_{f_{\mathbb{K}}}$
   has only one vertex, then
 Proposition 3.5 and Lemma 3.7 are true,  even if $\Delta$ is not a simplex. 
\end{Rem}

\begin{figure}[h]
\includegraphics[angle=0,width=0.3\textwidth]{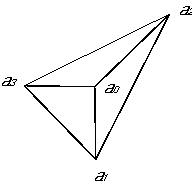}
\includegraphics[angle=0,width=0.3\textwidth]{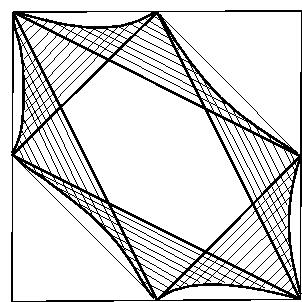}
\includegraphics[angle=0,width=0.3\textwidth]{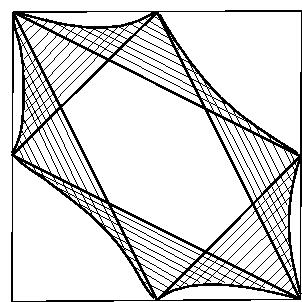}
\caption{Example of extra-pieces in dimension two for some choice of the coefficients, and the adjacent element of the subdivision}\label{c}
\end{figure}
\begin{figure}[h]
\includegraphics[angle=0,width=0.3\textwidth]{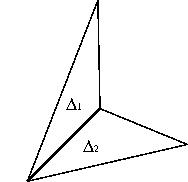}
\includegraphics[angle=0,width=0.6\textwidth]{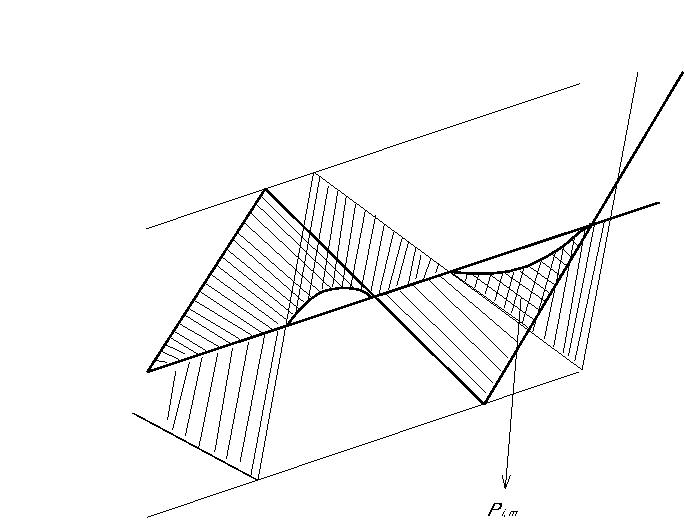}
\caption{Example of extra-pieces in dimension two, and the adjacent element of the subdivision}\label{c}
\end{figure}

\vspace{0.3cm}

\begin{Rem}
\begin{itemize}
\item[(i)]\,
The number of connected component of $\Arg
  (W(V_{f_{\mathbb{K},\, std}}))$, when we remove the real points, 
is $2^n-2$
and the volume of any component is $\frac{n-1}{n}\pi^n$,
\item[(ii)]\, if we denote by $\tilde{\mathscr{P}}_l$ the lift set in
  $\mathbb{R}^n$ of $\mathscr{P}_l$, then any component of
  $\tilde{\mathscr{P}}_l$ is a polyhedron (triangle for $n=2$ and not
  convex for $n>2$) with vertices in $\{ (k_1\pi ,\ldots ,k_n\pi )
  \}_{k_i\in\mathbb{Z}}$,
\item[(iii)]\, if we assume that $u$ can have  negative values then
  we have:
\begin{enumerate}
\item if $0<u\leq 1$ then we can choose the coefficients such that the argument of $W(V_{g_{\mathbb{K},\, u}})$ is
   constant and the tropical hypersurface $\Log\circ
  W(V_{g_{\mathbb{K},\, u}})$ vary,
\item if $u<0$ then the argument of $W(V_{g_{\mathbb{K},\, u}})$ varies
  and the tropical hypersurface $\Log\circ
  W(V_{g_{\mathbb{K},\, u}})$ is constant,
\item the set of arguments of $W(V_{f_{\mathbb{K}}})$ is called the
coamoeba of the complex tropical hypersurface $W(V_{f_{\mathbb{K}}})$
and is in the same time the limit of the coamoebas of some sequence of $J_t$-holomorphic hypersurfaces.
We describe the two last
  points with more details in the forthcoming papers \cite{N1-07} and \cite{N2-07}.
\end{enumerate}
\end{itemize}
\end{Rem}
We draw in figure 4 the set of argument  (known as the  coamoeba, for
more detail see \cite{N1-07}) of the curve in
$(\mathbb{C}^*)^2$ defined by the polynomial $f_1(z,w) = w^3z^2+wz^3+1$
where the matrix  ${}^tL_1^{-1}$ is equal to $\frac{1}{7}\left(
  \begin{array}{cc} 3&-1\\ -2&3\end{array}\right) $
and in figure 5 the coamoeba of the curve in $(\mathbb{C}^*)^2$
defined by the polynomial $f_2(z,w) = w^2z^2+z+w$
where the matrix  ${}^tL_2^{-1}$ is equal to $\frac{1}{3}\left(
  \begin{array}{cc} 1&1\\ -2&1\end{array}\right) $

\begin{figure}[h]
\includegraphics[angle=0,width=0.5\textwidth]{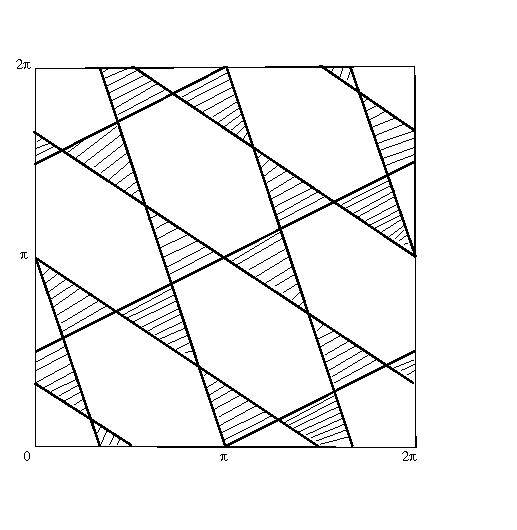}
\caption{The image of the curve defined by the polynomial $f_1(z,w)=wz^3+z^2w^3+1$ under the argument map $\Arg$}
\label{c}
\end{figure}

\begin{figure}[h]
\includegraphics[angle=0,width=0.4\textwidth]{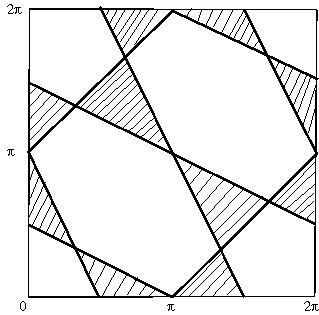}
\caption{The image of the curve defined by the polynomial $f(z,w)=z+w+z^2w^2$  under the argument map $\Arg$}
\label{c}
\end{figure}

\begin{figure}[h]
\includegraphics[angle=0,width=0.4\textwidth]{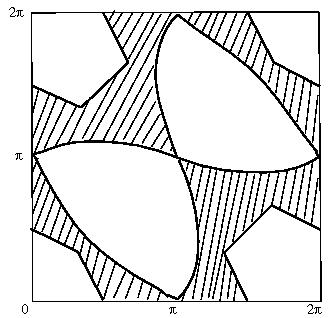}
\caption{The image of the curve defined by the polynomial $f(z,w)=z+w+\frac{zw}{2}+z^2w^2$ under the argument map $\Arg$}
\label{c}
\end{figure}

\section{Viro's patchworking principle.}

Let $\Delta$ be a convex integer polytope and $\tau =\cup _{v=1}^l
\Delta_v$ a convex integer subdivision of $\Delta$ (we can see
Viro's theory in \cite{V-90} for more details of this definition and
generally on the patchwork principle). This means that there
exists a convex piecewise affine linear map $\nu : \Delta
\longrightarrow \mathbb{R}$ so that:
\begin{itemize}
\item[(i)]\, $\nu_{\mid \Delta_v}$ is affine linear for each $v$,
\item[(ii)]\, if $\nu_{\mid U}$ is affine linear for some open set $U\subset
  \Delta$, then there  exists $v$ such that $U\subset \Delta_v$.
\end{itemize}
Let $\tilde{\Delta}$ be the extended polyhedral of $\Delta$ associated
to $\nu$, that is the convex hull of the set $\{  (\alpha ,u)\in
\Delta\times \mathbb{R} \mid \, u\geq \nu (\alpha ) \}$. For any
$\Delta_v\in \tau$, let $\lambda (x) = <x, a_v>+b_v$ be the affine linear
map defined on $\Delta$ such that $\lambda_{\mid \Delta_v} = \nu_{\mid
  \Delta_v}$ where $<,>$ is the scalar product in $\mathbb{R}^n$, \,
$a_v=(a_{v,\, 1},\ldots ,a_{v,\, n})\in\mathbb{R}^n$ and $b_v$ is a real number.
We put $\nu ' = \nu -\lambda$ and we define the family of polynomials
 $\{ f_t '\}_{t\in ]0,\, \frac{1}{e}]}$ by:
$$
f_t '(z) =\sum_{\alpha\in A}\xi_{\alpha}t^{\nu '(\alpha )}z^{\alpha}
$$
where $\xi_{\alpha}\in \mathbb{C}$. Then we have:
\begin{eqnarray}
f_t '(z) &=&t^{-b_v}\sum_{\alpha\in A}\xi_{\alpha}t^{\nu (\alpha
  )}(z_1t^{-a_{v,\, 1}})^{\alpha_1}\ldots (z_nt^{-a_{v,\,
  n}})^{\alpha_n}\nonumber \\
&=&t^{-b_v}f_t\circ \Phi^{-1}_{\Delta_v,\, t}(z)\nonumber
\end{eqnarray}
 where $f_t$  is the polynomial defined by:
$$
f_t (z) =\sum_{\alpha\in A}\xi_{\alpha}t^{\nu (\alpha )}z^{\alpha}
$$
 and $\Phi_{\Delta_v,\, t}$ is the self diffeomorphism of
 $(\mathbb{C}^*)^n$ defined by:
\[
\begin{array}{ccccl}
\Phi_{\Delta_v,\, t}&:&(\mathbb{C}^*)^n&\longrightarrow&(\mathbb{C}^*)^n\\
&&(z_1,\ldots ,z_n)&\longmapsto&(z_1t^{a_{v,\, 1}},\ldots ,z_nt^{a_{v,\,
  n}} ).
\end{array}
\]
This means that the polynomials $f_t '$ and $f_t\circ
\Phi^{-1}_{\Delta_v,\, t}$ defines the same hypersurface. So we have:
$$
V_{f_t '} = V_{f_t\circ
\Phi^{-1}_{\Delta_v,\, t}} = \Phi_{\Delta_v,\, t}(V_{f_t})
$$
Let $\Gamma_t$ be the
spine of the amoeba $\mathscr{A}_{H_t(V_{f_t})}$ where $H_t$
denotes the self diffeomorphism of $(\mathbb{C}^*)^n$ defined by
$H_h$ with $h=-\frac{1}{\log t}$ and $\Log_t = \Log \circ H_t$. 
Let $U(v)$ be a small ball in $\mathbb{R}^n$ with center the
vertex of $\Gamma_t$ dual to $\Delta_v$, $f_t^{\Delta_v}$ be  the truncation
  of $f_t$ to $\Delta_v$,
 and $V_{\infty ,\,\Delta_v}$ is the complex tropical hypersurface
    with tropical coefficients of index $\alpha\in
    \Delta_v$ (i.e., $V_{\infty ,\,\Delta_v} = \lim_{t\rightarrow 0}H_t(V_{f_t^{\Delta_v}})$).
Using Kapranov's  theorem (see \cite{K-00}),  we obtain the
 following proposition (called  tropical localization by Mikhalkin, 
see \cite{M2-04}):

\begin{Pro} For any $\varepsilon >0$ there exist $t_0$ such that if
  $t\leq t_0$ then the image under $\Phi_{\Delta_v,\, t}\circ H_t^{-1}$   of $H_t(V_{f_t})\cap \Log^{-1}(U(v))$ 
 is contained in the
    $\varepsilon$-neighborhood of the image under $\Phi_{\Delta_v,\, t}\circ H_t^{-1}$  
    of the complex tropical hypersurface
    $V_{\infty ,\,\Delta_v}$ ,
    with respect to the product metric
    in $(\mathbb{C}^*)^n\approx\mathbb{R}^n\times (S^1)^n$.
\end{Pro}

\begin{prooof} By decomposition of $f_t '$, we have:
$$
f_t '(z) =t^{-b_v}\sum_{\alpha\in \Delta_v\cap A}\xi_{\alpha}t^{\nu
  (\alpha ) -<\alpha ,a_v>}z^{\alpha} \,\,  +\,\, \sum_{\alpha\in
  A\setminus \Delta_v}  \xi_{\alpha}t^{\nu
  (\alpha ) -<\alpha ,a_v>-b_v}z^{\alpha} \qquad\qquad\qquad (5)
$$
On the other hand we have the following commutative diagram:
\begin{equation}
\xymatrix{
(\mathbb{C}^*)^n\ar[rr]^{\Phi_{\Delta_v,t}}\ar[d]_{\Log_t}&&
(\mathbb{C}^*)^n\ar[d]^{\Log_t}\cr
\mathbb{R}^n\ar[rr]^{\phi_{\Delta_v}}&&\mathbb{R}^n
} \qquad\qquad\qquad\qquad (6)     \nonumber
\end{equation}
such that if $v=(a_{v,\, 1},\ldots ,a_{v,\, n})\in \mathbb{R}^n$ is
the vertex of the tropical hypersurface $\Gamma$ dual to the
element $\Delta_v$ of the subdivision $\tau$,    then
$\phi_{\Delta_v}(x_1,\ldots , x_n)=(x_1-a_{v,\, 1},\ldots
,x_n-a_{v,\, n})$. Let $U(v)$  be a small open ball in
$\mathbb{R}^n$ centered at $v$.

\noindent Assume that $\Log_t(z)\in \phi_{\Delta_v}(U(v))$ and $z$ is no
singular in $V_{f_t}$. Then the second sum in $(5)$ converges to zero
when $t$ tends to zero, because by the choice of $z$ and $U(v)$, the
tropical monomials in $f_{trop,\, t}'$ corresponding to lattice points
of $\Delta_v$ dominates the monomials corresponding to lattice points
of $A\setminus \Delta_v$. But the first sum in $(5)$ is just a
polynomial defining the hypersurface $ \Phi_{\Delta_v,\,
    t} (V_{f_t^{\Delta_v}})$.

\noindent By the commutativity of the last diagram, if we take
$z\in V_{f_t'}$ such that $\Log_t(z)\in \phi_{\Delta_v}(U(v))$
then $\Log_t\circ \Phi_{\Delta_v,t}^{-1}(z)\in U(v)$ and hence
$H_t( \Phi_{\Delta_v,t}^{-1}(z))\in \Log^{-1}(U(v))$. So the image under 
$\Phi_{\Delta_v,\, t}\circ H_t^{-1}$ of
$H_t(V_{f_t})\cap \Log^{-1}(U(v))$ is contained in an
$\varepsilon$-neighborhood of  the image under  $\Phi_{\Delta_v,\, t}\circ H_t^{-1}$ of
  $H_t(V_{f_t^{\Delta_v}})$ for
sufficiently small $t$ and the proposition is done because
$V_{\infty ,\,\Delta_v}$ is the limit when $t$ tends to zero of
the sequence of $J_t$-holomorphic hypersurfaces
$H_t(V_{f_t^{\Delta_v}})$ (by taking a discreet sequence  $t_k$
converging to zero if necessary).
\end{prooof}

\section{Maximally sparse polynomials and proof of the main theorem}

From now we assume that the polynomial $f$ is maximally sparse
i.e. $\supp (f) = \Verte (\Delta_f)$. The family of polynomials $(4)$
can be considered as polynomial $f_{\mathbb{K}}$ with
coefficients in the non-Archimedean field $\mathbb{K}$ of Puiseux
series with coefficients in $\mathbb{C}$. So if we denote by
$V_{\mathbb{K}}$ the hypersurface in $(\mathbb{K}^*)^n$ defined by the
polynomial $f_{\mathbb{K}}$ and $-1/\log t$ the contraction of
$\mathbb{R}^n\longrightarrow\mathbb{R}^n$ defined by $(x_1,\ldots
,x_n)\longmapsto (-\frac{x_1}{\log t},\ldots ,-\frac{x_n}{\log t})$
for $t\in ]0,\frac{1}{e}]$, and $V_{f_t}$ the hypersurface defined by the the complex polynomial $f_t$,
 then we have the following theorem of
M. Passare and H. Rullg\aa rd in \cite{PR1-04} and G. Mikhalkin in \cite{M1-02}:

\begin{The} The non-Archimedean amoeba
  $\mathscr{A}_{V_{\mathbb{K}}} =\Gamma_{\infty} \subset\mathbb{R}^n$ of the
  hypersurface $V_{\mathbb{K}}\subset (\mathbb{K}^*)^n$ is the limit (
  with respect to the Hausdorff metric on compacts) of $(-1/\log
  t)(\mathscr{A}_{V_{f_t}})$ when $t$ tends to zero.
\end{The}

\vspace{0.1cm}

On one hand the non-Archimedean amoeba
  $\Gamma_{\infty}$ is the variety of the tropical
  polynomial $\displaystyle{\max_{\alpha\in \supp (f)}\{ -\nu (\alpha
  )+<\alpha , x>\}}$
and on the other hand the limit of $(-1/\log
  t)(\mathscr{A}_{V_{f_t}})$ is the limit of the spines $\Gamma_t$
 of the amoebas $\mathscr{A}_{H_t(V_{f_t})}$ of the $J_t$-holomorphic hypersurface $H_t(V_{f_t})$. Hence $\Gamma_{\infty}$ is
 solid ( because $\supp (f) = \Verte (\Delta_f)$ and any vertex of
 $\Delta_f$ corresponds to a complement component of the amoeba, see \cite{PR1-04})
 and the subdivision $\tau_{\infty}$ of $\Delta_f$ dual to
 $\Gamma_{\infty}$ has the following properties : by a small
 perturbation of the coefficient vector of $f$ if necessary, we can
 assume that the subdivision $\tau_{\infty}$ is a triangulation (this
 means that each element of $\tau_{\infty}$ is a simplex). Let $V_{f_{\vec a}}$ be the hypersurface defined by the coefficient vector $a =(a_1,\ldots , a_r)$, then
by the
 lower semi-continuity of the function $a\mapsto \sharp\{
 \mbox{component of $\mathbb{R}^n\setminus \mathscr{A}_{f_{\vec a}}$}\}$, if
 the coefficient vector $\vec b$ is  close enough to $\vec a$,  then the
 number of complement components of $\mathbb{R}^n\setminus
 \mathscr{A}_{f_{\vec b}}$ is greater or equal to the
 number of complement components of $\mathbb{R}^n\setminus
 \mathscr{A}_{f_{\vec a}}$ (see \cite{FPT-00}). Hence if we prove that $ \mathscr{A}_{f_{\vec b}}$ is
 solid then $ \mathscr{A}_{f_{\vec a}}$ is solid too. So we can suppose for
 our problem that $\tau_{\infty}$ is a triangulation.

\vspace{0.3cm}

Note that the vertices of any simplex of $\tau_{\infty}$ are contained in $\Verte (\Delta_f)$.
Let $\mathscr{L}\subset \Delta_f\cap\mathbb{Z}^n$ be the complement of
$\Verte (\Delta_f)$ in the image of the order mapping
i.e. $\mathscr{L} = \{  \alpha\in\Delta_f\cap\mathbb{Z}^n \mid \,
\, E_{\alpha}^c\ne \emptyset \, \mbox{and}\, a_{\alpha} =
0\}$.
Using the triangulation  $\tau_{V_f}$ dual to the spine of the
amoeba $\mathscr{A}_{V_f}$, we define a triangulation $\tau_{V_f,\,
\mathscr{L}}$ of $\Delta_f$  satisfying the following properties (see Appendix. C for more details and notations) :
\begin{itemize}
\item[(i)]\, $\tau_{V_f,\, \mathscr{L}} = \tau_{\infty}$,
\item[(ii)]\, there is a deformation $\{ H_t(V_{f_t})\}_{t\in ]0,\, \frac{1}{e}]}$
 of the hypersurface $V_f$ so that for any $t\in ]0,\, \frac{1}{e}]$  the spines
 $\Gamma_t$ of the amoebas $\mathscr{A}_{H_t(V_{f_t})}$ have the same combinatorial
 type than $\Gamma_f$, the spine of the amoeba of the initial hypersurface $V_f$; and
$\displaystyle{\lim_{t\rightarrow 0} \Gamma_t = \Gamma_{\infty}}$. This means that the triangulation of $\Delta_f$ dual
 to $\Gamma_t$ is the same for each $t\in ]0,\, \frac{1}{e}]$ and  the triangulation
 of $\Delta_f$ dual to $\Gamma_{\infty}$  is $\tau_{V_f,\, \mathscr{L}}$.
\end{itemize}

\vspace{0.2cm}

Our main aim in this section is to show that, the amoebas of the $J_t$-holomorphic hypersurfaces $H_t(V_{f_t})$'s cannot develop a new complement components in $\mathbb{R}^n$ other than those of $\Gamma_{\infty}$. More precisely, we show that the spines of those amoebas have the same combinatorial type that the one of $\Gamma_{\infty}$.

\begin{The} Let $\mathscr{S} = \{ t\in ]0,\, \frac{1}{e}] \mid \,
  \mbox{the\, amoeba}\,\, \mathscr{A}_{H_t(V_{f_t})}\,\,\mbox{is\,
    solid} \}$. Then $\mathscr{S}$ is a nonempty, closed and open
  subset  of $]0,\, \frac{1}{e}]$, and hence equal to $]0,\, \frac{1}{e}]$.
\end{The}

\begin{Lem} For a sufficiently small $t$ the amoebas
  $\mathscr{A}_{H_t(V_{f_t})}$ are solid. In particular $\mathscr{S}$
  is nonempty.
\end{Lem}

\begin{prooof}Assume that there exists an infinite sequence $\{t_m\}$
  which tends to zero and  such that the amoebas
  $\mathscr{A}_{H_{t_m}(V_{f_{t_m}})}$ are not solid, and the order of the complement component of the amoebas is $\beta\in \Delta_{i}$, where $\Delta_{i}$ is an element of the subdivision $\tau_{\infty}$ of $\Delta$.

This means that
  there exists a sequence of parallel hyperplanes
  $\mathscr{P}_m\subset \mathbb{R}^n\times\mathbb{R}$ dual to the
  vertex $\beta$,
  such that $\mathscr{P}_m\cap \cup_{j=1}^l\mathscr{P}_{\alpha_j}$ is equal to the union of compact
  polyhedrons in $\Gamma_{t_m}$ the spine of
  $\mathscr{A}_{H_{t_m}(V_{f_{t_m}})}$, and $\mathscr{P}_{\alpha_j}$
  are the hyperplanes of $\mathbb{R}^n\times\mathbb{R}$ dual to the
  vertices $\alpha_j$'s of $\Delta_{i}$ if
  $\beta\in Int(\Delta_{i})$. The hyperplanes $\mathscr{P}_{\alpha_j}$
  of $\mathbb{R}^n\times\mathbb{R}$   are the hyperplanes  dual to the
  vertices $\alpha_j$'s of
  $\Delta_{i}\cup \Delta_{l}$ if
  $\beta\in \Delta_{i}\cap
  \Delta_{l}$. Hence
  for any $n\geq 2$, the
  deformation can have $n$ possibilities:
\begin{itemize}
\item[(i)]\, if the order $\beta$ of the new complement component
  (i.e. of order not in $\Verte (\Delta_f)$) is in the interior of
  some $\Delta_{i}$ we have the first possibility,
\item[(ii)]\, and if the order of the new complement component is
  contained in a face
   of $\Delta_{i}$ we have
  $n-1$ possibilities, one possibility for any positive dimension of the faces
  of $\Delta_{i}$.
\end{itemize}
We can see the two possibilities when $n=2$ in figures 7 and 8.

\begin{figure}[h]
\includegraphics[angle=0,width=0.25\textwidth]{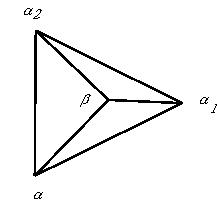}\qquad\includegraphics[angle=0,width=0.25\textwidth]{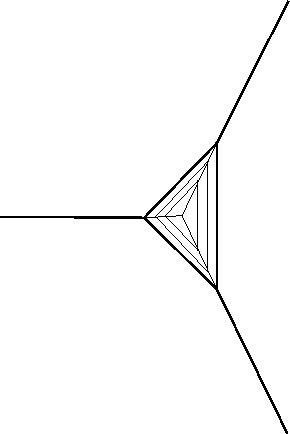}
\caption{$n=2$ and $\beta\in Int (\Delta_{\mathscr{C}_{\alpha}})$}
\label{c}
\end{figure}

\begin{figure}[h]
\includegraphics[angle=0,width=0.25\textwidth]{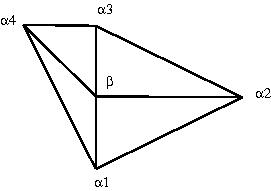}\qquad\includegraphics[angle=0,width=0.25\textwidth]{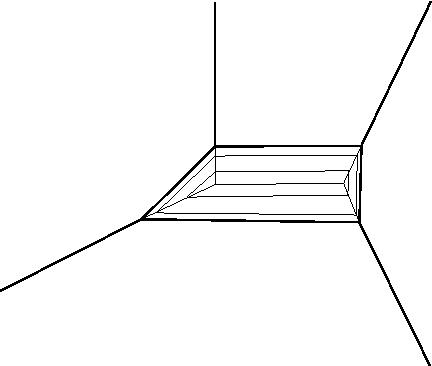}
\caption{$n=2$ and $\beta\in \partial \Delta_{\mathscr{C}_{\alpha}}$}
\label{c}
\end{figure}

\begin{figure}[h]
\includegraphics[angle=0,width=0.25\textwidth]{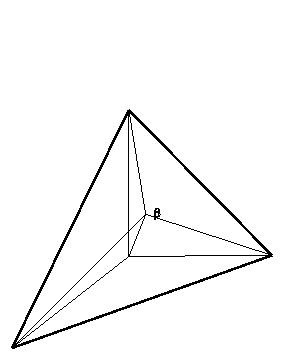}\quad\includegraphics[angle=0,width=0.25\textwidth]{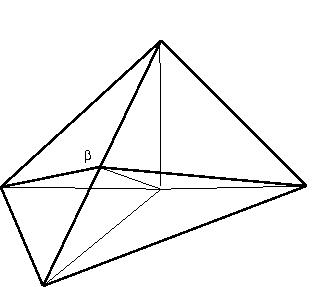}
\quad\includegraphics[angle=0,width=0.25\textwidth]{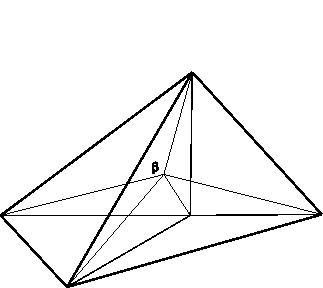}
\caption{The three possibilities for $n=3$}
\label{c}
\end{figure}

Let $\displaystyle{\mathscr{P}_0=\lim_{m\rightarrow \infty}\mathscr{P}_m}$ which
is a hyperplane parallel to the $\mathscr{P}_m$'s. If $\beta\in
Int(\Delta_{i})$, then   either  the $\mathscr{P}_m$'s  passes
through the lifting in $\mathbb{R}^n\times\mathbb{R}$ of the vertex $v$ of $\Gamma_{\infty}$ dual to
$\Delta_{i}$ or 
they
go to infinity . But if $\beta\in \partial
(\Delta_{i})$, then $\mathscr{P}_0$ is an
hyperplane parallel to the $\mathscr{P}_m$'s and
containing the dual of the sub-simplex of $\partial
(\Delta_{i})$ in which $\beta$ lies , or
$\mathscr{P}_0$ is parallel to the $\mathscr{P}_m$'s and goes to
infinity.
We can treat only the first case; the others can be given in
the same way if we restrict ourself to the sub-simplex in the
boundary of $\Delta_{i}$ containing $\beta$.

\vspace{0.1cm}

In other words, if we denote by $\nu_{\infty}(\beta )$ the limit
of $\nu_m(\beta )$ when $m$ tends to infinity ( i.e. $t_m$ tends
to zero), with $\nu_m$  the Passare-Rullg\aa rd function
corresponding to the spine of the amoeba $\mathscr{A}_{f_{t_m}}$,
  we have a priori two possibilities :
\[
\lim_{m\rightarrow \infty} \nu_m(\beta ) = \left\{ \begin{array}{c}
-\infty   \qquad \qquad \qquad \qquad  (a)\\
\mbox{or}\\
<\beta ,a_v> + b_v\qquad\qquad  (b)
\end{array}
\right. \]
where $<\beta ,a_v> + b_v$ is the equation of the hyperplane in
$\mathbb{R}^n\times\mathbb{R}$
passing through the points of coordinates $(\alpha_j, \nu (\alpha_j))$
with $\alpha_j\in\Verte (\Delta_{i})$.

Let $g_t^{(m)}(z) = f_t(z) + t^{\nu_m (\beta)} (t-t_m)z^{\beta}$ and
$V_{\infty ,\, m}$ be the complex tropical hypersurface equal to the
limit of $H_t(V_{g_t^{(m)}})$ when $t$ tends to zero. Using  theorem 5.1, it's clear that
$\Log (V_{\infty ,\, m}) = \Gamma_{\infty ,\, m}$ is the tropical
hypersurface equal to the spine of the amoeba
$\mathscr{A}_{H_{t_m}(V_{f_{t_m}})}$  and $\displaystyle{\lim_{m\rightarrow
  \infty} \Gamma_{\infty ,\, m} = \Gamma_{\infty}}$.

\vspace{0.1cm}

\noindent Let us show that case (a) cannot occur.
If $\displaystyle{\lim_{m\rightarrow
    \infty}\nu_{m}(\beta ) =-\infty }$, then  $\displaystyle{\lim_{m\rightarrow
    \infty}c_{\beta ,\, m} =\infty }$ (the constants defined by Passare and Rullg\aa rd, see Section 2) and we have two cases:
\begin{itemize}
\item[(i)]\, either the hyperplane of $\mathbb{R}^{n+1}$ corresponding
  to $\beta$
intersect no other
  hyperplane corresponding to $\alpha\in \Verte (\Delta )$, and then
  the monomial of index $\beta$ is omitted from   $g_t^{(m)}$ for each $m$;
and hence $\Log^{-1}(v)\cap V_{\infty ,\, f} =\Log^{-1}(v)\cap  V_{\infty ,\,
    0}$,  with $v$  the vertex in  the non-Archimedean  amoeba $\displaystyle{\Gamma_{\infty} = \lim_{m\rightarrow \infty}\mathscr{A}_{H_{t_m}(V_{g_t^{(m )}})}}$
 dual to ${\Delta}_{i}$,
$V_{\infty ,\, f}$ is the complex tropical hypersurface
    equal to the limit of $H_t(V_{f_t})$ when $t$ tends to zero, and
    $V_{\infty ,\, 0}$ is the limit of $V_{\infty ,\, m}$ when $m$ tends
     to infinity (with respect to the Hausdorff metric).
\item[(ii)]\, or the amoeba $\Gamma_{\infty} $ 
has a complement
  component of order $\beta$ such that $c_{\beta ,\, \infty}= + \infty$ and the 
  coefficient $b_{\beta ,\, m}$ of $g_t^{(m)}$ of index $\beta$ evaluated at $0$ is
  not bounded i.e. tends to $\infty$. This case cannot occur
  because it contradicts the fact that the amoeba
  $\Gamma_{\infty}$ is solid and  nonempty.
\end{itemize}

\vspace{0.2cm}

\noindent In   case (b), by multiplying $f_t$ by $t^p$ such that
$p+b_v>0$ if necessary, we can assume that $b_v>0$. Recall that $\displaystyle{V_{\infty ,\, f} = \lim_{t\rightarrow 0}H_t(V_{f_t})}$,\,$\displaystyle{V_{\infty ,\, m} = \lim_{t\rightarrow 0}H_t(V_{g_t^{(m)}})}$  and $\displaystyle{V_{\infty ,\, 0} = \lim_{m\rightarrow \infty} V_{\infty ,\ m}}$.

\vspace{0.2cm}
  
\begin{Pro} In case (b) we have 
$V_{\infty ,\, f} = V_{\infty ,\,
    0}$, in particular, if $v$ is the vertex of $\Gamma_\infty$ dual to the simplex $\Delta_i$ containing $\beta$, then  $\Log^{-1}(v)\cap  V_{\infty ,\,
    0} = \Log^{-1}(v)\cap  V_{\infty ,\,
    f}$.
\end{Pro}

The problem is only on the 0-dimensional cell $v$ of the non-Archimedean amoeba $\Gamma_\infty$ dual to the simplex ${\Delta}_{i}$  of the triangulation and  containing $\beta$. We denote by $\delta_{l}^k$ the $k$-dimensional cells of $\Gamma_\infty$ containing $v$ as vertex with $k=1,\ldots ,n-1$, which are the  dual to the $(n-k)$-faces $F_l^{n-k}$ of ${\Delta}_{i}$ of positive dimension. If $\tau_m$ is the triangulation of $\Delta$ dual to $\Gamma_{\infty ,\, m}$, and $\Delta_l\subset{\Delta}_{i}$ are the elements of $\tau_m$ of maximal dimension, then we denote by $\delta_{l}^{(m),\, k}$ (resp. $\delta_{l}^k$) 
the $k$-cells of $\Gamma_{\infty ,\, m}$ (resp. of $\Gamma_\infty$) which are  dual to the $(n-k)$-faces of the proper faces $F_l^{n-1}$ of ${\Delta}_{i}$ for $k=1,\ldots , n-1$ (see M. Passare and H. Rullg\aa rd \cite{PR1-04} for more details).

\begin{Lem}For all $l$ and $k$ we have : 
$$\displaystyle{\lim_{m\rightarrow \infty}(V_{\infty ,\, m}\cap \Log^{-1}(Int (\delta_{l}^{(m),\, k}))) = V_{\infty ,\, f}\cap \Log^{-1}(Int (
\delta_{l}^k))},
$$ 
where $Int (\delta_{l}^{(m),\, k})$ is the interior of $\delta_{l}^{(m),\, k}$.
\end{Lem}

\begin{prooof}
Let $U(\delta_{l}^{(m),\, k})\subset \mathbb{R}^n$ be a small open neighborhood of the interior
of $\delta_{l}^{(m),\, k}$ satisfying the following properties: (i) its intersection with any other cells of $\Gamma_m$ of dimension less than $k$ is empty, and (ii) the limit of $U(\delta_{l}^{(m),\, k})$ when $m$ tends to infinity is $Int (\delta_{l}^{k})$ (see \cite{V-90} page 42 for more details). Firstly we know that 
$$
V_{\infty ,\, m}\cap \Log^{-1}(Int(\delta_{l}^{(m),\, k})) = \lim_{t\rightarrow 0} H_t(V_{g_t^{(m)}})\cap \Log^{-1}(U(\delta_{l}^{(m),\, k})).
$$
This means that for any $m$ and any $\varepsilon  > 0$, there exists $T_m<\frac{1}{e}$ such that if $t\leq T_m$ then $H_t(V_{g_t^{(m)}})\cap \Log^{-1}(U(\delta_{l}^{(m),\, k}))$ is contained in an $\varepsilon$-neighborhood of $V_{\infty ,\, m}\cap \Log^{-1}(Int(\delta_{l}^{(m),\, k}))$.  We look now at $V_{f_t}$ as the end of a path $\gamma$ of hypersurfaces in $(\mathbb{C}^*)^n$ starting at $V_{g_t^{(m)}}$, where the parameter of the path $\gamma$ is the valuation of the coefficient of index $\beta$. This means that the coefficients of index different than $\beta$ are independent of the parameter. By the continuity of roots property (see for example \cite{B-71}) we have for any $\varepsilon >0$ there exists $\eta >0$ such that if $\mid \nu_m (b_\beta ) - \nu_{m'} (b_\beta ) \mid < \eta$ then the Hausdorff distance between $H_t(V_{g_t^{(m)}})\cap \Log^{-1}(U(\delta_{l}^{(m),\, k}))$ and $H_t(V_{f_t})\cap \Log^{-1}(U(\delta_{l}^{(m'),\, k}))$ is less than $\varepsilon$ for a sufficiently small $t$ (we can assume that there exists  a very large $m'$ such that $t = t_{m'}$). On the other hand $\displaystyle{V_{\infty ,\, f}\cap \Log^{-1}(Int (\delta_{l}^k)) = \lim_{m\rightarrow \infty} H_{t_m}(V_{f_{t_m}}) \cap \Log^{-1}(U(\delta_{l}^{(m),\, k}))}$.   By the triangular inequality of the Hausdorff distance, we obtain that for any $\varepsilon ' >0 $, there exists $m_1$ such that if $m > m_1$ then we have
$$
d_{\mathcal{H}}(V_{\infty ,\, m}\cap \Log^{-1}(Int(\delta_{l}^{(m),\, k})); V_{\infty ,\, f}\cap \Log^{-1}(Int (\delta_{l}^k))) < \varepsilon ',
$$
 and the lemma is done. Proposition 5.4 is an immediate consequence of Lemma 5.5, and in particular, if $v$ is the vertex of $\Gamma_\infty$ dual to the simplex ${\Delta}_{i}$ containing $\beta$, then  $\Log^{-1}(v)\cap  V_{\infty ,\,
    0} = \Log^{-1}(v)\cap  V_{\infty ,\,
    f}$.

\begin{figure}[h]
\includegraphics[angle=0,width=0.7
\textwidth]{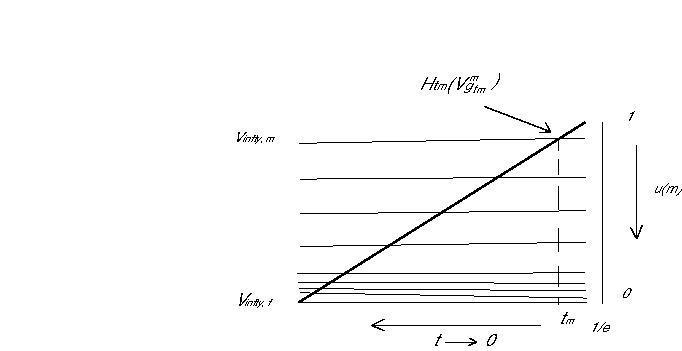}
\caption{}
\label{c}
\end{figure}

\end{prooof}

\vspace{0.1cm}

\noindent {\em End of proof of Lemma 5.3.} The triangulation of $\Delta_f$, dual to the spine $\Gamma_{t_m}$ of the amoeba of the hypersurface $H_{t_m}(V_{g^{(m)}_t})$  is unchanged (i.e. the same for each $m$), because the $\Gamma_{t_m}$'s have the same combinatorial type. 
Indeed the set of slopes of faces of a tropical hypersurface  is a
finite set of rational numbers.
Hence, by a restriction to a subsequence of $\{ t_m\}$ if necessary, we may assume that the valuation of the coefficients $b_{\beta ,\, m}$'s take their values in the interval $[-\nu_{m_0}(\beta ), -\nu_{\infty}(\beta )]$ for some $m_0$, so that there exists a strictly decreasing function $u: \{m\in\mathbb{N}/ m>m_0\}\rightarrow [0,1]$ with $\displaystyle{\lim_{m\rightarrow \infty}u(m)=0}$ and $u(m_0)=1$. So the valuation $\val (a_{\beta ,\, m})=- \nu_m(\beta )$ can be written as follow : $(1-u(m))(-\nu_{\infty}(\beta )) + u(m)(-\nu_{m_0}( \beta ))$. So the  case $(b)$  satisfy the hypothesis $(i)$ and $(ii)$ which are the assumptions of Lemma 3.7.

\noindent By Proposition 4.1, if $U(v)\subset \mathbb{R}^n$ is a small ball  centered on $v$, then the  tropical localization says that 
for any $\varepsilon >0$ there exist $t_0$ such that if
$t\leq t_0$ then $\Arg (H_t(V_{f_t})\cap \Log^{-1}(U(v)))$ is contained in the
$\varepsilon$-neighborhood $\mathscr{N}_{\varepsilon}(\Arg (V_{\infty ,\,\Delta_{i}}))$ of the set of arguments of the complex tropical hypersurface $V_{\infty ,\, \Delta_{i}}$ (because the transformations $H_t$ and $\Phi_{{\Delta}_{i},\, t }$ conserve the arguments); this means that $\Arg (V_{\infty ,\, f}\cap \Log^{-1}(v))$ is contained in an  $\varepsilon$-neighborhood of  $\Arg (V_{\infty ,\,\Delta_{i}})$, where $\displaystyle{V_{\infty ,\,\Delta_{i}} =\lim_{t\rightarrow 0} H_t(V_{f_t^{\Delta_{i}}})}$.

\begin{figure}[h]
\includegraphics[angle=0,width=0.4\textwidth]{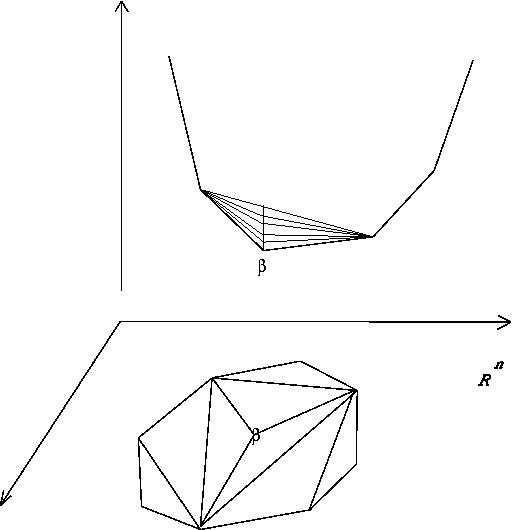}
\caption{The deformation given by the valuation}
\label{c}
 \end{figure}

\noindent On the other hand $\Arg ( H_t(V_{g_t^{(m,\, \Delta_{i})}})\cap \Log^{-1}(U(v)))$
 is contained in the
    $\varepsilon$-neighborhood $\mathscr{N}_{\varepsilon}(\Arg ((V_{\infty ,\, \Delta_{i}}')))$ of the set of arguments of the complex tropical hypersurface $V_{\infty ,\, \Delta_{i}}'$ which is the limit of the complex tropical hypersurfaces $V_{\infty ,\,g_m^{\Delta_{i}}}$ when $m$ tends to $\infty$, and $\displaystyle{V_{\infty ,\,g_m^{\Delta_{i}}} = \lim_{t\rightarrow 0}H_t(V_{g_t^{(m,\, \Delta_{i})}})}$ (here $g_t^{ (m,\, \Delta_{i})}$ is the truncation of $g_t^{(m)}$ to $\Delta_{i}$), and then $\Arg (V_{\infty ,\, 0}\cap\Log^{-1}(v))$ is contained in an $\varepsilon$-neighborhood of  $\Arg ( V_{\infty ,\, \Delta_{i}}')$.

\noindent If $\mathcal{P}_j$ is an extra-piece in $\Arg ( V_{\infty ,\, \Delta_{i}}')$ (see Lemma 3.7), then we claim that $\mathcal{P}_j\cap \Arg (V_{\infty ,\, 0}\cap\Log^{-1}(v))$ has a no vanishing volume. Indeed, assume that $\Vol (\mathcal{P}_j\cap \Arg (V_{\infty ,\, 0}\cap\Log^{-1}(v))) = 0$, hence there exists  an external hyperplane $H_{ij}$ for $\Arg ( V_{\infty ,\, \Delta_{i}}')$ which is not external for $\Arg ( V_{\infty ,\, 0})$, such that $H_{ij}\cap \partial (\overline{\Arg (V_{\infty ,\, 0})})$ is of dimension $n-1$. But this situation cannot occur, because the hyperplane $H_{ij}$ is not external for $\Arg ( V_{\infty ,\, 0})$.

\noindent The set of arguments $\Arg (V_{\infty ,\, f}\cap \Log^{-1}(v))$ is contained in an  $\varepsilon$-neighborhood of  $\Arg (V_{\infty ,\,\Delta_{\mathscr{C}_{\alpha}}})$, and its intersection with 
$\mathcal{P}_j$ is empty (recall that the polynomial $f$ is maximally sparse, see Lemma 3.7). Hence $V_{\infty ,\, f}$ and $V_{\infty ,\, 0}$ cannot be equal, and then  we have
a contradiction, which  means that such sequence of $t_m$'s cannot 
exist. Hence for sufficiently small $t$, the amoebas
$\mathscr{A}_{H_t(V_{f_t})}$ are solid and then the set $\mathscr{S}$ is nonempty.
\end{prooof}

\vspace{0.1cm}

\noindent The following Corollary is a consequence of the last construction and Proposition 3.5. 
\begin{Cor}Let $V_f\subset (\mathbb{C}^*)^n$ be an hypersurface
defined by a maximally sparse polynomial $f$ with Newton polytope
a simplex. Then the amoeba of $V_f$ is solid.
\end{Cor}

\begin{prooof} Assume instead that the amoeba $\mathscr{A}_f$  is not solid;
hence there exist $\beta\in \Delta\cap\mathbb{Z}^n$ which is the order
 of some complement component other than those  corresponding to
 $\Verte (\Delta )$. Let 
$g_{t}^{(m)}(z) = f_{t_m}
 (z) + e^{(1-u(m))s+u(m)\nu (\beta )} t^{(1-u(m))s+u(m)\nu (\beta )}(t_m - t)z^{\beta}$
  where $t_m$ is a sequence of real numbers which tends to zero, and $s = <\beta , a_v> + b_v$, where 
$y=<x , a_v> + b_v$ is the equation of the hyperplane
  in $\mathbb{R}^n\times \mathbb{R}$ containing the points of coordinates $(\alpha ,
  \nu (\alpha ))$ with $\alpha \in \Verte (\Delta_f)$, and the sequence $u(m)$ is the sequence defined above.
Using the above deformation and  applying Proposition 3.5, we obtain that the complex tropical hypersurfaces  $\Arg (V_{\infty ,\, f})$ and
 $\Arg (V_{\infty ,\, 0})$ are different, because even if $s\gg
<\beta , a_v> + b_v$ and tends to infinity, the set of arguments
$V_{\infty ,\, 0}$ contains extra-pieces
corresponding to the no vanishing coefficients which have no
contribution on the non-Archimedean amoeba.

\end{prooof}

\noindent {\it Proof of theorem}.5.2.   
By Lemma 5.3, the set
$\mathscr{S}$ is nonempty and it is obviously closed. Let $t_{max}$ be
the maximum of $t\in \mathscr{S}$. We claim  that $t_{max}$ is in  the interior of
the interval $]0; \frac{1}{e}]$ and then $\mathscr{S}$ is open. Indeed, assume on the contrary that there is an infinite  sequence $\{ t_m\}$ in $]t_{max}; \frac{1}{e}]$ such that $\lim_{m\rightarrow \infty}t_m = t_{max}$, and the amoebas of the 
hypersurfaces $V_m = \{ z\in (\mathbb{C}^*)^n\mid \, f_{t_m}(z)=0 \}$ are not solid. We know that the amoebas of the hypersurfaces defined by the truncated polynomials $f_t^{\Delta_i}$ are solid, because the $\Delta_i$'s are a simplexes, and the set of its arguments contains no extra-pieces. For sufficiently large $m$, let we assume that  $f_{t_m}$ develop just one complement component of order $\beta$, and $\Delta_i$ is the element of $\tau$ containing $\beta$. Let $g_t^{(m)}$ be the family of polynomials defined by   $g_t^{(m)}(z) = f_t(z) + t^{\nu (\beta )}(t-t_m) z^{\beta}$.
Let us denote by $g_t^{(m;\, \Delta_i)}$ the truncation of $g_t^{(m)}$ to $\Delta_i$, and $V_{(m;\, \Delta_i)}$ the complex tropical hypersurface which is the limit of $H_t(V_{g_t^{(m;\, \Delta_i)}})$ when $t$ tends to zero (with respect to the Hausdorff metric on compacts). For any $m$, using the same reasoning as in Lemma 3.7, we show that the  complex tropical hypersurface  $V_{(\infty ;\, \delta_i)}$ equal to the limit of  $V_{(m;\, \Delta_i)}$ when $m$ tends to the infinity, contains extra-pieces with no vanishing volume. This means that the set of arguments of $H_{t_{max}}(V_{g_{t_{max}}^{(\infty ;\, \Delta_i)}})$ contains some extra-pieces. Contradiction because $H_{t_{max}}(V_{g_{t_{max}}^{(\infty ;\, \Delta_i)}})$ is equal to $H_{t_{max}}(V_{f_{t_{max}}^{\Delta_i}})$ (see Proposition 5.4). Hence there is no sequence $\{ t_m \}$ such that
 $\displaystyle{t_m\rightarrow t_{max}}$ and the amoebas
 $\mathscr{A}_{H_{t_m}(V_{f_{t_m}})}$ are not solid, and then $\mathscr{S}$ is open.

\begin{Cor} Let $V$ be an algebraic hypersurface in $(\mathbb{C}^*)^n$
  defined by a maximally sparse polynomial $f$ with amoeba
  $\mathscr{A}_{V}$ and  Newton polytope
  $\Delta$. Then the number of  components of $
  \mathbb{R}^n\setminus \mathscr{A}_{V}$ is equal to the number of
  vertices of $\Delta$.
\end{Cor}


\section{Appendix}

\vspace{0.3cm}

\centerline{\bf A: Proposition 3.5 in the Case  $n=1$}

\vspace{0.3cm}

Let us prove in this Appendix the Lemma 3.6 and proposition 3.5 in  the  case   $n=1$.
Let $\Delta = [0, k]$,\,
      $A=\{\beta_1,\,\ldots ,\beta_s\}\subset \, ]0,k[\, \cap\,
      \mathbb{Z}$. Firstly, we can remark that if $f_{\mathbb{K}}(z)
      =a_kz^k+a_0$, then we can   assume that the coefficient
      $a_k$  is equal to  one and the valuation
      of  the coefficient  $a_0$ is zero.
\begin{figure}[h]
\includegraphics[angle=0,width=0.6\textwidth]{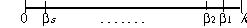}
\label{c}
\end{figure}
 Indeed, the first assertion
      is obvious and for the second we consider the translation  $\Phi_{a_0}$ of
      $\mathbb{K}^*$ defined by
\[
\begin{array}{ccccl}
\Phi_{a_0}&:&\mathbb{K}^*&\longrightarrow&\mathbb{K}^*\\
&&z&\longmapsto&t^{-\frac{\val (a_0)}{k}}z
\end{array}
\]
and we obtain:
\[
\begin{array}{ccccl}
 f_{\mathbb{K}}\circ\Phi_{a_0}(z) &=& t^{-\val (a_0)}z^k + a_0\\
                                  &=& t^{-\val (a_0)} (z^k+a_0')\\
                                  &=&t^{-\val (a_0)}f_{\mathbb{K}}'(z)
\end{array}
\]
with $\val ( a_0') = 0$ and
$\Phi_{a_0}^{-1}(V_{f_{\mathbb{K}}}) =  V_{f_{\mathbb{K}}'}$. So
we have:
\[
\begin{array}{ccccl}
W(V_{f_{\mathbb{K}}})&=&e^{\frac{-\val (a_0)+i\arg (a_0)}{k}}.\{ e^{i\frac{(2l+1)\pi}{k}}\}_{l=0}^{k-1}   \\
\end{array}
\]
Let $g_{\mathbb{K},\, u}(z) = z^k+ a_{\beta_1,\, u}z^{\beta_1}+\cdots
+a_0$ such that the initial part of $a_{\beta_1,\, u}$ is $t^{-\val
  (a_{\beta_1,\, 1})u} e^{i\arg (a_{\beta_1,\, u})}$ and set

$g_{\mathbb{K},\, u}^{[\beta_1,\, k]}= z^k+ a_{\beta_1,\, u}z^{\beta_1}$ (the truncation of
$g_{\mathbb{K},\, u}$ to $[\beta_1,\, k]$).
Therefore
$W(V_{g_{\mathbb{K},\, u}^{[\beta_1,\, k]}}) = e^{\frac{-\val (a_{\beta_1,\, 1})u +
    i\arg (a_{\beta_1,\, 1})}{k-\beta_1}}.\{
e^{i\frac{(2l+1)\pi}{k-\beta_1}}\}_{l=0}^{k-\beta_1-1}$ and for each
$j=1,\ldots ,s$ we have in a similar way:
$$
W(V_{g_{\mathbb{K},\, u}^{[\beta_j,\, \beta_{j-1}]}}) = e^{A_j+iB_j}.\{
e^{i\frac{(2l+1)\pi}{\beta_{j-1}-\beta_j}}\}_{l=0}^{\beta_{j-1}-\beta_j-1}
$$
where $A_j = \frac{-\val (a_{\beta_{j,\, 1}})u+\val (a_{\beta_{j-1,\,
      1}})u}{\beta_{j-1}-\beta_j}$ and
      $B_j =  \frac{\arg (a_{\beta_{j,\, 1}})-\arg (a_{\beta_{j-1,\,
      1}})}{\beta_{j-1}-\beta_j}$. So we obtain:
$$
W(V_{g_{\mathbb{K},\, 0}}) = \bigcup_{j=1}^{s+1}e^{A_j+iB_j}.\{
e^{i\frac{(2l+1)\pi}{\beta_{j-1}-\beta_j}}\}_{l=0}^{\beta_{j-1}-\beta_j-1}
$$
with $a_{\beta_{0,\, 1}}= a_k,\,\,a_{\beta_{s+1,\, 1}} = a_0$ and $\{
\beta_{j}\}_{j=1}^s = A$. The first part of the lemma is done if we
put $v_j=B_j$ for $j=1,\ldots , s+1$ and $L_j : x \mapsto
(\beta_{j-1}-\beta_j)x$.\\

If $n=1$ then  we can see that $W(V_{g_{\mathbb{K},\,
  0}})= W(V_{f_{\mathbb{K}}})$ if and only if $\beta_{j-1}-\beta_j$
is a constant $r$ and $k=r(s+1)$. So it
suffices to prove the assertion for the polynomials
$f_{\mathbb{K}}°(z)=z^{s+1}+a_0$ and $g_{\mathbb{K},\,u}°(z)
=z^{s+1}+a_0+\sum_{j=1}^sa_{\beta_{j,\, u}}z^{j}$.
By an easy computation we have $a_{\beta_1,\, 0}\ne
0$, so  in this case the
roots $r_j$ of the polynomial $g_{\mathbb{K},\,u}°$ cannot have the
same absolute value,
 because if it is the case then their sum is zero (because of the
 condition on their arguments). But  in
the case when the absolute value of the roots $r_j$ are not the
same then the amoeba of the limit of the $g_{\mathbb{K},\,u}°$,
when $u$ tends to zero, have at least two points, which contradict
the fact that the limit of the  non-Archimedean amoeba has only
one point.

\vspace{0.3cm}

\centerline{\bf  B:  The Set of Arguments of the Standard Hyperplane}

\vspace{0.3cm}

let $\mathscr{P}_{sdt}$ be the hyperplane in $(\mathbb{C}^*)^n$
defined by the polynomial $f(z_1,\ldots ,z_n) = 1+\sum_{i=1}^n
z_i$ with Newton polytope the standard simplex. If $(S^1)_l^{n-1}$
is the $(n-1)$-torus in $(S^1)^n$ defined by: $\{ (x_1,\ldots
,x_{l-1},e^{i\pi},x_{l+1},\ldots ,x_n)\in (S^1)^n  \}$ for
$l=1,\ldots ,n$, then the lift in $\mathbb{R}^n$ of the union
$\displaystyle{\cup_{l=1}^n(S^1)_l^{n-1}}$ divide a fundamental
domain of the torus into $2^n$ parts $\{ \tau_s \}_{s=1}^{2^n}$.
Let $\tau_1$ be the $n$-cube in $\mathbb{R}^n$ of vertices
$(v_1,\ldots ,v_{n-1},\pi )$, $(v_1,\ldots ,v_{n-1},2\pi )$ with
$v_j\in\{ 0,\pi \}$, and $\mathscr{C}_1$ is the cone of vertex
$v_0=(0,\ldots ,0,\pi )$ and
 base the $(n-1)$-cube of vertices $(v_1,\ldots ,v_{n-1},2\pi )$ with  $v_j\in\{ 0,\pi \}$.

\begin{Lem}
The image under the argument map $\Arg$ of the hyperplane $\mathscr{P}_{std}\subset (\mathbb{C}^*)^n $
defined by the equation $1+\sum_{i=1}^n z_i=0$
 is the union of the $2^n-2$ polyhedron $\mathscr{D}_s=\tau_s\setminus \mathscr{C}_s$ (not convex for $n>2$)   such that :
\begin{itemize}
\item[(a)]\, The $\tau_s$ are different than the two following cubes:
(i) $\tau_0$ of vertices  $(v_1,\ldots ,v_{n-1},\pi )$, $(v_1,\ldots ,v_{n-1},0 )$ with $v_j\in\{ 0,\pi \}$, and (ii) $\tau_{\pi}$ of vertices  $(v_1,\ldots ,v_{n-1},\pi )$ and \\
 $(v_1,\ldots ,v_{n-1},2\pi )$ with $v_j\in\{ \pi ,2\pi \}$,
\item[(b)]\, the polyhedrons $\mathscr{D}_s$ viewed as subset of the
 torus (by projection), are two by two attached by $2^n-1$ real points
 of $(S^1)^n$, and the complement in $(S^1)^n$ of the closure of  their union is a
 connected and convex polyhedron,
\item[(c)]\, the volume of the coamoeba of $\mathscr{P}_{std}$ is
equal to $\frac{(n-1)(2^n-2)}{n}\pi^n$ (with respect to the  flat metric of the torus).
\end{itemize}
\end{Lem}

\begin{prooof} If $n=2$, then we have $z_2 =- (1+z_1)$, so \,
$\arg (z_2) = \pi + \arg (1+z_1) \mod (2\pi )$. Hence if $\arg
 (z_1) =\alpha < \pi$ and its module varies between zero and the infinity,
 then $\arg (z_2)$ varies  between $\pi$ and $\pi +\alpha$. By switching the
  role of the variable $z_1$ and $z_2$, the lemma is done (see figure 8).
  For $n>2$, we use induction on $n$ and the fact that  $\displaystyle{z_n = -
   (1+ \sum_{j=1}^{n-1} z_j)}$. Put $\alpha_j = \arg (z_j)$, so if $0\leq \alpha_j\leq \pi$
   for $j=1,\ldots n-1$ then $\alpha_n\in ]\pi , m_n[$ where $m_n = \max_{1\leq j\leq n-1}\{ \alpha_j\}$
   and then we have one of the sets $\mathscr{D}_s$. By changing the position of the arguments of the $z_j$'s
   we obtain the $2^n - 2$ sets. Convexity is local property, which is the case in our situation,
   so the second statement of the lemma is obvious. For the third affirmation of the
   lemma, it suffice to compute the volume of the cone's $\mathscr{C}_s$ which is equal
   to $\pi^n/n$ and then the volume of any $\mathscr{D}_s$ is $\pi^n - \pi^n/n = \frac{(n-1)}{n}\pi^n$.  I leave the details to the reader.

\begin{figure}[h]
\includegraphics[angle=0,width=0.4\textwidth]{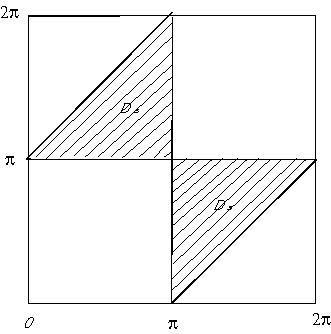}
\caption{The left picture represent one $\mathscr{D}_s=
\tau_s\setminus \mathscr{C}_s$ and the right one represent all
the argument of a line which has two $\mathscr{D}_s$'s subsets.}
\label{c}
\end{figure}

\begin{figure}[h]
\includegraphics[angle=0,width=0.6\textwidth]{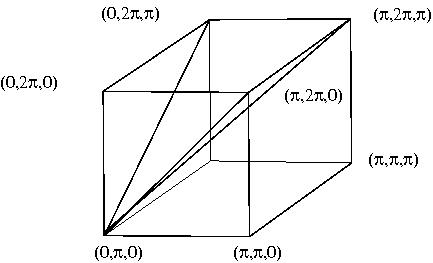}
\caption{If $n=3$, we represent here one $\mathscr{D}_s=\tau_s\setminus \mathscr{C}_s$
 where $\tau_s$ is all the cube, and  $\mathscr{C}_s$ is the cone of base the square in the
 top and vertex the point $(0,\pi ,0)$.}
\label{c}
\end{figure}

\end{prooof}

\begin{Rem} We can see, using the result of section 3, that if $g\in \mathbb{K}[z_1,\ldots ,z_n]$ is a maximally sparse polynomial with Newton polytope a simplex $\Delta$, and  $V_g$ is the hypersurface in $(\mathbb{K}^*)^n$ defined by $g$, then
the volume of the set of arguments of $W(V_g)$  is equal to the volume of the set of argument of the standard hyperplane $\mathscr{P}_{std}$  in $(\mathbb{C}^*)^n$ defined by the polynomial $\displaystyle{f(z_1,\ldots ,z_n) = 1+\sum_{i=1}^n z_i}$ (the torus is equipped  with the flat metric). Indeed, $\Arg (W(V_g)) = {}^tL_{\Delta}^{-1}\Arg (\mathscr{P}_{std})$ (viewed in some fundamental domain in $\mathbb{R}^n$ the universal covering of the torus), where $L_{\Delta}$ is the linear part of the affine linear surjection $\rho :\Delta_{std}\rightarrow \Delta$. Firstly we have $\Vol (\mathscr{D}_s) = \frac{\Vol (\mathscr{D}_s')}{\det (L_{\Delta})}$, where $\mathscr{D}_s'$ are the polyhedrons in the torus $(S^1)^n$  corresponding to $W(V_g)$, and for any $\mathscr{D}_s$ it corresponds $\det (L_{\Delta})$ times $\mathscr{D}_s'$. Hence we have $\Vol (\Arg (W(V_g))) = \Vol (\cup \mathscr{D}_s') = \det (L_{\Delta}) \left(\frac{\cup \mathscr{D}_s}{\det (L_{\Delta})}\right) = \Vol (\Arg (\mathscr{P}_{std}))$.

\end{Rem}


\vspace{0.3cm}

\centerline{\bf C: Construction of the Triangulation $\tau_{V_f,\, \mathscr{L}}$}

\vspace{0.3cm}

We use the notation of Section 5 where  $\mathscr{L}\subset \Delta_f\cap\mathbb{Z}^n$ denotes  the complement of
$\Verte (\Delta_f)$ in the image of the order mapping
i.e. $\mathscr{L} = \{  \alpha\in\Delta_f\cap\mathbb{Z}^n \mid \,
\mbox{with}\, E_{\alpha}^c\ne \emptyset \, \mbox{and}\, a_{\alpha} =
0\}$;  the polynomial $f$ is assumed maximally sparse.  Using the triangulation $\tau_{V_f}$ dual to the spine of the
amoeba $\mathscr{A}_{V_f}$, we define a new triangulation $\tau_{V_f,\,
\mathscr{L}}$ of $\Delta_f$ as follow:

\begin{itemize}
\item[{\it Step}.1]\,  Let $\alpha_1\in \mathscr{L}$ and denote by
  $\Delta_{\mathscr{C}_{\alpha_1}^j}$ the following subsets of
  $\tau_{V_f}$ :
\begin{enumerate}
\item \, $\Delta_{\mathscr{C}_{\alpha_1}^1} = \cup_i\Delta_i$ is a
  convex subset of $\Delta_f$ containing $\alpha_1$ where $\Delta_i$'s
  are element of $\tau_{V_f}$,
\item \, $\Verte (\Delta_{\mathscr{C}_{\alpha_1}^1})\subseteq \Verte
  (\Delta_f)$,
\item \, for any proper face $F_{\Delta_{\mathscr{C}_{\alpha_1}^1}}$
  of $\Delta_{\mathscr{C}_{\alpha_1}^1}$ we have $\partial
  F_{\Delta_{\mathscr{C}_{\alpha_1}^1}}\subset\cup_j\partial
  F^j_{\Delta_f}$ where $\cup_j F^j_{\Delta_f} = \partial \Delta_f$
\item \, $\Delta_{\mathscr{C}_{\alpha_1}^1}$ is of minimal volume with
  properties (1), (2) and (3).
\end{enumerate}
If there exist $r$ subsets of $\tau_{V_f}$ satisfying (1), (2), (3)
and (4)
with $r\geq 3$ and  the interior of
$\Delta_{\mathscr{C}_{\alpha_1}^i}\cap\Delta_{\mathscr{C}_{\alpha_1}^j}$
is empty, then we associate  to $\alpha_1$ the union
$\displaystyle{\cup_{j=1}^r\Delta_{\mathscr{C}_{\alpha_1}^j}}$. We can
remark that this case can occur only if $\alpha_1$ is contained in the
boundary of a proper face of $\Delta_f$.\\
If $\alpha_1$ is in the interior of a proper face of
$\Delta_{\mathscr{C}_{\alpha_1}^1}$, then either there is many other
$\Delta_{\mathscr{C}_{\alpha_1}^j}$'s satisfying (1),(2),(3) and (4),
in this case we associate to $\alpha_1$ their union, or there is only
one $\Delta_{\mathscr{C}_{\alpha_1}^1}$  satisfying (1),(2),(3) and
(4) and then we take the connected component of $\Delta_f\setminus
\Delta_{\mathscr{C}_{\alpha_1}^1}$ containing $\alpha_1$ and we repeat
the same operation on this component. So we obtain
$\Delta_{\mathscr{C}_{\alpha_1}^j}$'s simplex  with $\Vol
(\Delta_{\mathscr{C}_{\alpha_1}^j})\geq \Vol
(\Delta_{\mathscr{C}_{\alpha_1}^1})$.
As in the second case we associate to $\alpha_1$ the union of
$\Delta_{\mathscr{C}_{\alpha_1}^1}$ with the
$\Delta_{\mathscr{C}_{\alpha_1}^j}$'s.
If $\alpha_1\in Int (\Delta_{\mathscr{C}_{\alpha_1}^1})$ then there is
a unique subset $\Delta_{\mathscr{C}_{\alpha_1}^1}$ of $\tau_{V_f}$
satisfying (1), (2), (3) and (4).

\begin{figure}[h]
\includegraphics[angle=0,width=0.5\textwidth]{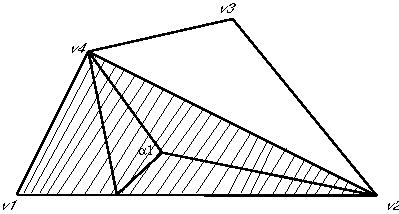}
\caption{The dashed triangle represent one $\Delta_{\mathscr{C}_{\alpha_1}^1}$}
\label{c}
\end{figure}

\begin{figure}[h]
\includegraphics[angle=0,width=0.3\textwidth]{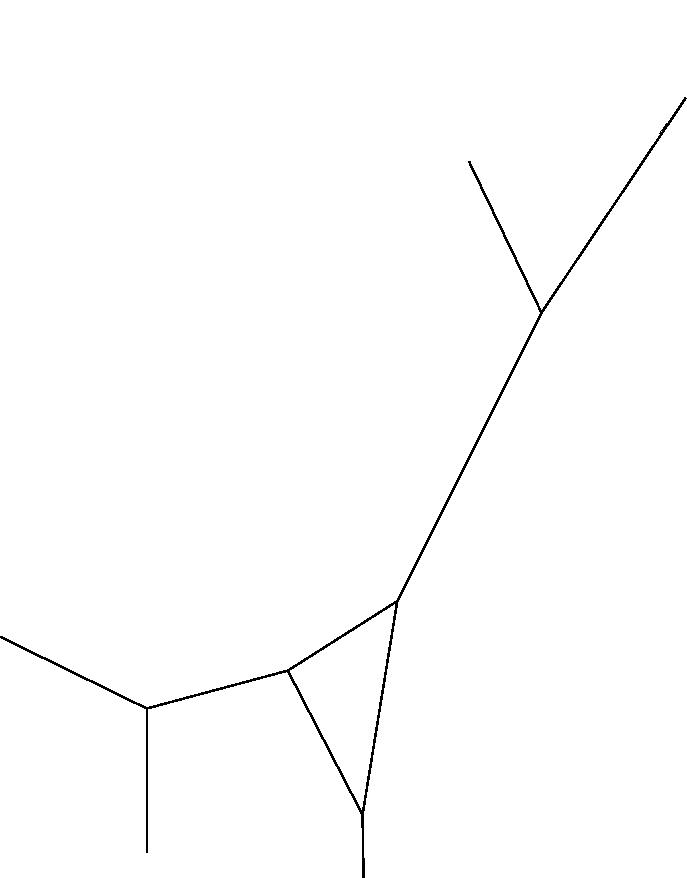}
\includegraphics[angle=0,width=0.3\textwidth]{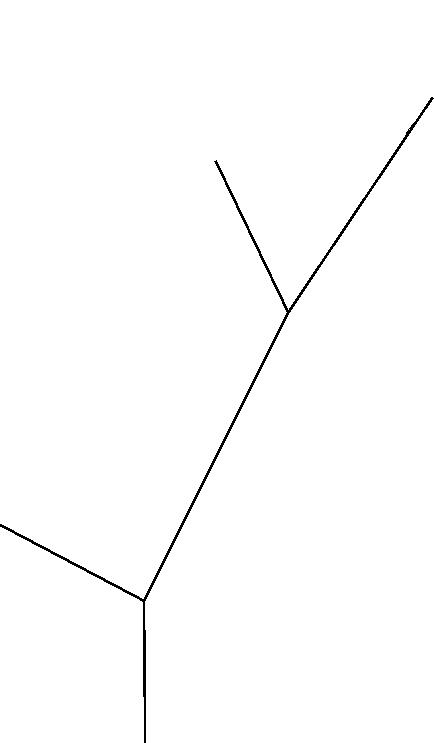}
\caption{The spines $\Gamma_t$ have the same combinatorial type
 of the left picture and $\Gamma_0$ has the combinatorial type of the right one}
\label{c}
\end{figure}

\item[{\it Step}.2]\, Let $K_{\alpha_1} = \Delta_f$ and $K_{\alpha_2}$
  be the connected component of $\Delta_f\setminus
  (\cup_j\Delta_{\mathscr{C}_{\alpha_1}^j})$ containing $\alpha_2$
and we repeat the same operation for $\alpha_2$. By this process we
obtain a new subdivision $\tau_{V_f,\, \mathscr{L}}$ of $\Delta_f$ such
that $\tau_{V_f,\, \mathscr{L}} = \mathscr{P}\cup \mathscr{R}$ where
$\mathscr{P} =
\cup_{i=1}^s\cup_{j=1}^{r_i}\Delta_{\mathscr{C}^j_{\alpha_i}}$ and
$\mathscr{R}$ is the union of element in $\tau_{V_f}$ not in $\mathscr{P}$.
\end{itemize}

Let $\Gamma_{\mathscr{C}^j_{\alpha_i},\, t}$ be the tropical
hypersurface dual to ${\tau_{V_f}}_{\mid
  \Delta_{\mathscr{C}^j_{\alpha_i}}}$ such that all its vertices are
in $\Verte (\Gamma_t)$ where $\Gamma_t$ is the spine of the amoeba
$\mathscr{A}_{H_t(V_{f_t})}$. We denotes by
$\mathscr{C}^j_{{\alpha_i},\, t}$ the set of its bounded polyhedrons
of dimension $n-1$. If $\alpha\in \mathscr{L}$ we denotes, in the sequel,
 by $\mathscr{C}_{\alpha ,\, t}$ one of the $\mathscr{C}^j_{{\alpha},\, t}$'s.
\begin{Def}
 We said that a family of a polyhedron $P_t\subset \mathbb{R}^n$
vanishes or collapses if we have the following:
\begin{itemize}
\item[(i)]\, for each $t>0$, the polyhedrons  $P_t$ are
  homothetic,
\item[(ii)]\, their  volume  tends to zero
  when $t\longrightarrow 0$.
\end{itemize}
\end{Def}
This means that the dual of $P_t$ is constant for each $t>0$ and $P_t$
collapses to some point.

\begin{Lem}
Let $\alpha\in \mathscr{L}$. Then the set $\mathscr{C}_{\alpha ,\, t}$ vanishes.
\end{Lem}

\begin{prooof} If the amoebas $\mathscr{A}_{H_t(V_{f_t})}$ converge
  (with respect to the Hausdorff metric in the compact subsets of
  $\mathbb{R}^n$) to a tropical hypersurface $\Gamma_{\infty}$, then also,
  their spines converge to $\Gamma_{\infty}$. In particular the number
  of polyhedrons in $\Gamma_{\infty}$ of maximal dimension
  (i.e. $n-1$) is not greater than the number of polyhedrons of
  $\Gamma_t$ of maximal dimension. This means that some polyhedrons
  $P_t\subset \Gamma_t$ converge to a  parallel polyhedron $P\subset
  \Gamma_{\infty}$
(because the set of slopes of faces of a tropical hypersurface  is a   
finite set of rational numbers)
and some other one vanished.\\
Let $\alpha\in\mathscr{L}$ and
$\Gamma_{\mathscr{C}_{\alpha} ,\, t}$ be the dual of
${\tau_{V_f}}_{\mid \Delta_{\mathscr{C}_{\alpha}}}$ for
$\alpha\in\mathscr{L}$. By definition, the set of vertices of
$\Gamma_{\mathscr{C}_{\alpha} ,\, t}$ is contained in $\Verte
(\Gamma_t)$, on the other hand, if  $K$ is a compact in $\mathbb{R}^n$
containing all vertices of $\Gamma_{\mathscr{C}_{\alpha} ,\, t}$, then
$\Gamma_t\cap\Gamma_{\mathscr{C}_{\alpha} ,\, t}$ converges to the
intersection of $K$ with the tropical hypersurface
$\Gamma_{\Delta_{\mathscr{C}_{\alpha}}}$ dual to
$\Delta_{\mathscr{C}_{\alpha}}$ (we mean here the dual to the
subdivision of the polytope $\Delta_{\mathscr{C}_{\alpha}}$ with only
one element). Hence $\mathscr{C}_{\alpha ,\, t}$ vanishes, because
$\Gamma_{\Delta_{\mathscr{C}_{\alpha}}}$ has only one vertex, and then
it has no compact subpolyhedrons other than its vertex.
\end{prooof}

\begin{Rem}
We can remark that
 the subdivision $\tau_{V_f,\, \mathscr{L}}$ is a convex
  subdivision of $\Delta_f$ and defined by the Passare-Rullg\aa rd's function
$\nu$ restricted to the set of vertices of $\Delta$. So $\tau_{V_f,\,
    \mathscr{L}}$ is dual to $\Gamma_{\infty}$ and hence $\tau_{V_f,\,
    \mathscr{L}} = \tau_{\infty}$.
\end{Rem}

\vspace{0.3cm}

\centerline{\bf D: Maximally sparse is an optimal condition}

\vspace{0.3cm}

We give in this Appendix an example of a curve $V_f\subset (\mathbb{C}^*)^2$ defined by a polynomial $f$  of Newton polygon $\Delta$ such that its support contains other elements than those of the  vertices of $\Delta$, and the number of the complement components of the amoeba $\mathscr{A}_f$ is strictly greater than the cardinality of the support of $f$.
Let $V_f$ be the curve in $(\mathbb{C}^*)^2$ defined by the polynomial  $f(z,w) = -zw^2 + z^3w -7zw+6w +z$. We can see, by computation that the points $(0,\,0),\, (\log (2),0)$ and $(\log (3),\, 0)$ are contained in the amoeba $\mathscr{A}_f$ of the curve $V_f$ and the two points $(\frac{\log (2)}{2},\, 0)$ and $(\frac{\log (3)}{2},\, 0)$ are contained in two compact different complement components of $\mathscr{A}_f$. Hence the number of complement components  in $\mathbb{R}^2$ of $\mathscr{A}_f$is equal to $6$ which is strictly greater than the number of monomials of $f$.

\begin{figure}[h]
\includegraphics[angle=0,width=0.2\textwidth]{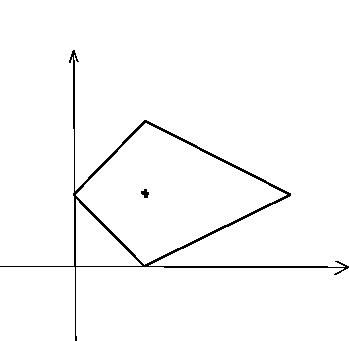}\quad
\includegraphics[angle=0,width=0.2\textwidth]{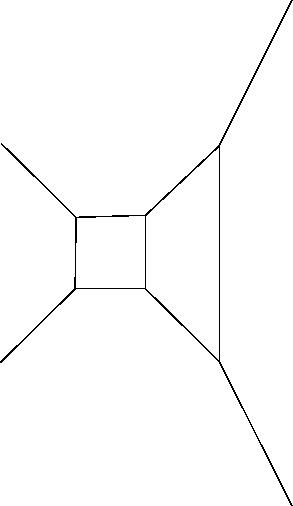}\quad
\includegraphics[angle=0,width=0.2\textwidth]{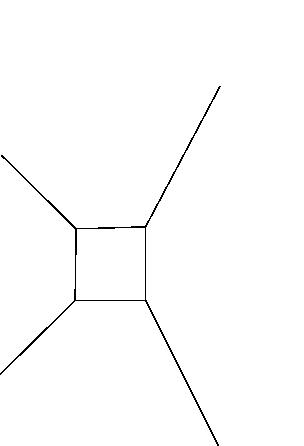}
\caption{The Newton polygon of $f$,the spine of $\mathscr{A}_f$, and the non-Archimedean amoeba $\Gamma_{\infty}$}
\label{c}
\end{figure}


\begin{thebibliography}{9}

\bibitem[B-71]{B-71}{\sc G. M. Bergman}, {\em The logarithmic limit-set of an algebraic variety},
Trans. Amer. Math. Soc.   {\bf 157}, (1971), 459-469


\bibitem[FPT-00]{FPT-00}{\sc M. Forsberg, M; Passare
and A. Tsikh}, {\em Laurent determinants and arrangements of hyperplane amoebas},
Advances in Math.  {\bf 151}, (2000), 45-70.


\bibitem[GKZ-94]{GKZ-94}{\sc I. M. Gelfand, M.
M. Kapranov and A. V. Zelevinski}, {\em Discriminants, resultants and
multidimensional determinants},
Birkh{\"a}user Boston 1994.

\bibitem[K-00]{K-00}{\sc M. M. Kapranov}, {\em Amoebas over
non-Archimedean fields}, Preprint 2000.

\bibitem[M1-02]{M1-02}{\sc G. Mikhalkin}, {\em  Decomposition into
pairs-of-pants for complex algebraic hypersurfaces},Topology {\bf 43}, (2004), 1035-1065.


\bibitem[M2-04]{M2-04}{\sc G. Mikhalkin}, {\em  Enumerative Tropical
Algebraic Geometry In $\mathbb{R}^2$}, J. Amer. Math. Soc. {\bf 18}, (2005), 313-377.

\bibitem[N1-07]{N1-07}{\sc M. Nisse}, {\em  Coamoebas of complex algebraic
    hypersurfaces}, in preparation.

\bibitem[N2-07]{N2-07}{\sc M. Nisse}, {\em  Coamoebas of curves in $(\mathbb{C}^*)^2$, Dimer Models on Torus and Planar Quivers}, in preparation.


\bibitem[PR1-04]{PR1-04}{Passare and Rullg\aa rd}{\sc M. Passare and H.
Rullg\aa rd}, {\em Amoebas, Monge-Amp{\`e}re measures, and triangulations
of the Newton polytope}, Duke Math. J. {\bf 121}, (2004), 481-507.


\bibitem[PR2-01]{PR2-01}{\sc M. Passare and H.
Rullg\aa rd}, {\em Multiple Laurent series and polynomial amoebas},
pp.123-130 in: Actes des rencontres d'analyse complexe, Atlantique,
{\'E}ditions de l'actualit{\'e} scientifique, Poitou-Charentes 2001.

\bibitem[R-01]{R-01}{\sc H. Rullg\aa rd}, {\em Polynomial
amoebas and convexity}, Research Reports In Mathematics Number 8,2001,
Department Of Mathematics Stockholm University.

\bibitem[V-90]{V-90}{\sc O. Viro}, {\em  Patchworking real
 algebraic varieties}, preprint: http://www.math.uu.se/ oleg; Arxiv: AG/0611382

\end{thebibliography}
\end{document}